\documentclass[11pt,reqno]{amsart}
\allowdisplaybreaks

\usepackage[T1]{fontenc}
\usepackage{lmodern}
\usepackage[a4paper,margin=1.15in]{geometry}
\usepackage{amsmath,amssymb,amsthm}
\usepackage{microtype}
\usepackage{etoolbox}
\makeatletter
\patchcmd{\@setaddresses}{\/:\space}{\/\space}{}{}
\patchcmd{\@setaddresses}{\/:\space}{\/\space}{}{}
\patchcmd{\@setaddresses}{\/:\space}{\/\space}{}{}
\makeatother
\usepackage{url}
\usepackage{hyperref}
\hypersetup{
    colorlinks=true,
    citecolor=blue,
    linkcolor=blue,
    urlcolor=blue,
    pdftitle={On two Romanoff type problems of Erdos},
    pdfauthor={Yuchen Ding},
    pdfkeywords={Romanoff theorem, primes, squarefree integers, additive complements, real powers, uniform distribution, covering congruences, Lucas numbers}
}

\numberwithin{equation}{section}

\newtheorem{theorem}{Theorem}[section]
\newtheorem{proposition}[theorem]{Proposition}
\newtheorem{lemma}[theorem]{Lemma}
\newtheorem{conjecture}[theorem]{Conjecture}

\title[On two Romanoff type problems of Erd\H{o}s]
{On two Romanoff type problems of Erd\H{o}s}

\author[Y. Ding]{Yuchen Ding}
\address{(Yuchen Ding) School of Mathematical Sciences, Yangzhou University, Yangzhou 225002, People's Republic of China\\ HUN-REN Alfr\'ed R\'enyi Institute of Mathematics, Budapest, Pf. 127, H-1364 Hungary}
\email{ycding@yzu.edu.cn}
\thanks{ORCID \href{https://orcid.org/0000-0001-7016-309X}{0000-0001-7016-309X}}

\keywords{Romanoff theorem, primes, squarefree integers, additive complements, real powers, uniform distribution, covering congruences, Lucas numbers}
\subjclass[2020]{11P32, 11B13, 11B34, 11K60, 11N36}

\begin{document}

    \begin{abstract}
        Let $\mathcal{P}$ be the set of primes and, for $a>1$, put
        $$
        \mathcal{S}_a=\{p+\lfloor a^k\rfloor\mid p\in\mathcal{P},\ k\in\mathbb{N}\}.
        $$
        Erd\H{o}s recorded a question of Kalm\'ar asking whether $\mathcal{S}_a$ has positive lower asymptotic density for every real $a>1$. We prove this for almost every real $a>1$. In fact,
        $$
        \liminf_{N\to\infty}\frac{|\mathcal{S}_a\cap[1,N]|}{N}
        \ge \frac{1}{\log a+9C_0/\pi^2},
        $$
        where $C_0$ is an absolute constant. The dependence on $a$ has the correct order $1/\log a$ as $a\to\infty$. The proof is based on residue distribution and weighted pair correlation estimates for the sequence $\lfloor a^k\rfloor$.

        At the same time, for almost every $a>1$ and every $\eta>0$, the number of positive integers $n\le x$ outside $\mathcal{S}_a$ is at least $x^{1-\eta}$ for all sufficiently large $x$. For almost every $a>1$, a main step is to construct, for all sufficiently large $K$ running through powers of two, congruences covering $\lfloor a^m\rfloor$, $1\le m\le K$, with a squarefree modulus $Q_K$ satisfying $\log Q_K=o(K)$. For the golden ratio $\varphi=(1+\sqrt5)/2$, we prove the explicit bounds
        $$
        0<\underline d(\mathcal{S}_\varphi)\le \overline d(\mathcal{S}_\varphi)\le \frac{1937}{1938}.
        $$

        We then turn to another problem of Erd\H{o}s, which asks whether every sufficiently large odd integer is the sum of a squarefree integer and a power of two. Let $\mathcal{Q}$ denote the positive squarefree integers. Replacing $2^m$ by $\lfloor a^m\rfloor$, we prove that for almost every $a>1$ and every $\varepsilon>0$, the exceptional sets for
        $$
        \mathcal{Q}+\{\lfloor a^m\rfloor\mid m\ge1\}
        \quad\text{and}\quad
        \mathcal{Q}+\{\lfloor a^m\rfloor+\lfloor a^\ell\rfloor\mid 1\le m<\ell\}
        $$
        are respectively
        $$
        O_{a,\varepsilon}\!\left(\frac{x(\log\log x)^{1+\varepsilon}}{\sqrt{\log x}}\right)
        \quad\text{and}\quad
        O_{a,\varepsilon}\!\left(\frac{x(\log\log x)^{1+\varepsilon}}{\log x}\right).
        $$
    \end{abstract}

    \maketitle

    \section{Introduction}

    Romanoff \cite{Romanoff} proved that, for every integer $a>1$, the set
    $$
    \{p+a^k\mid p\in\mathcal{P},\ k\in\mathbb{N}\}
    $$
    has positive lower asymptotic density. Erd\H{o}s and Tur\'an \cite{ErdosTuran1935a,ErdosTuran1935b} and Erd\H{o}s \cite{Erdos1951} later gave simpler proofs. The theorem led to a number of problems concerning sums of primes and sparse sequences. The case of powers of two has received particular attention. Quantitative estimates were obtained by Chen and Sun \cite{ChenSun}, Habsieger and Roblot \cite{HabsiegerRoblot}, L\"u \cite{Lv}, Pintz \cite{Pintz}, Habsieger and Sivak-Fischler \cite{HabsiegerSivak}, and Elsholtz and Schlage-Puchta \cite{ElsholtzSchlage}. Numerical work on the density of integers of the form $p+2^k$ was carried out by Romani \cite{RomaniThesis,Romani} and by Del Corso, Del Corso, Dvornicich and Romani \cite{DelCorso}.

    Erd\H{o}s also studied obstructions to such representations. Van der Corput \cite{Corput} proved that a positive proportion of the odd integers are not of the form $p+2^k$. Erd\H{o}s \cite{Erdos1950} constructed an odd arithmetic progression with the same property by means of covering congruences. Further explicit constructions were given by Habsieger and Roblot \cite{HabsiegerRoblot} and Chen, Dai and Li \cite{ChenDaiLi}. Crocker \cite{Crocker} and Pan \cite{Pan} studied the corresponding problem with two powers of two.

    In 1961 Erd\H{o}s \cite[p. 230, Problem (14)]{Erdos1961} recorded a question of Kalm\'ar in which the integral powers are replaced by integer parts of real powers. For $y>1$, put
    $$
    \mathcal{S}_y=\{p+\lfloor y^k\rfloor\mid p\in\mathcal{P},\ k\in\mathbb{N}\}.
    $$
    Erd\H{o}s asked whether $\mathcal{S}_y$ has positive lower asymptotic density for every $y>1$. The problem also appears as Problem 244 in Bloom's collection \cite{Bloom244}. The choice $k\ge1$ instead of $k\ge0$ is immaterial for lower asymptotic density, since the omitted set $\mathcal{P}+1$ has density zero. The distinction between almost every $y$ and every $y$ appears to be substantial. Tao observed in a comment on Problem 244 \cite{Tao244} that, for a fixed base, the integers $\lfloor y^k\rfloor$ could in principle be strongly concentrated in one residue class modulo many moduli, which could obstruct a positive density argument of Romanoff type. He further noted that excluding such concentration appears to be related to further progress on Mahler's classical problem on the powers of $3/2$ \cite{Mahler}. The metric distribution estimates used here rule out this behavior for almost every $y$, but do not address an arbitrary fixed base. Some individual irrational bases can nevertheless be treated by recurrence methods. Dubickas \cite{Dubickas} proved, for example, that the integers of the form $p+\lfloor(2+\sqrt3)^n\rfloor$ have positive lower asymptotic density. Ballot and Luca \cite[Theorem~1]{BallotLuca} proved a general Romanoff type theorem for nonconstant nondegenerate integral linear recurrences with separable characteristic polynomial. Our first result answers Kalm\'ar's question for almost every real base and gives a quantitative lower bound with the correct order of dependence on the base.

    \begin{theorem}\label{thm-romanoff}
        For $y>1$, put
        $$
        \delta_y=\liminf_{N\to\infty}\frac{\left|\mathcal{S}_y\cap[1,N]\right|}{N}.
        $$
        Let $C_0$ be the absolute constant defined in \eqref{eq-sieve-constant}. Then, for almost every $y>1$,
        \begin{equation}\label{eq-explicit-density-bound}
            \delta_y\ge \frac{1}{\log y+9C_0/\pi^2}.
        \end{equation}
        In particular, $\mathcal{P}+\{\lfloor y^k\rfloor\mid k\ge1\}$ has positive lower asymptotic density for almost all $y>1$.
    \end{theorem}

    The dependence on $y$ in \eqref{eq-explicit-density-bound} has the correct order as $y\to\infty$. Indeed, the prime number theorem \cite[Ch.~14]{Apostol} gives
    $$
    \begin{aligned}
        \left|\mathcal{S}_y\cap[1,N]\right|
        \le \sum_{k\le (\log N)/(\log y)+O_y(1)}\pi(N)\le (1+o_y(1))\frac{N}{\log y}.
    \end{aligned}
    $$
    Thus a lower bound which depends only on $y$ cannot have order larger than $1/\log y$.

    There are also explicit irrational bases satisfying the qualitative conclusion of Theorem~\ref{thm-romanoff} for a simple reason. For $m\ge2$, put
    $$
    A_m=\left\lfloor\left(\frac32\right)^m\right\rfloor,
    \qquad
    y_m=A_m^{1/m}.
    $$
    Since $A_m/(3/2)^m\to1$, we have $y_m\to3/2$. Moreover, $1<y_m<3/2$, and the rational root theorem shows that $y_m$ is irrational. Since
    $$
    \mathcal{P}+\{A_m^j\mid j\ge1\}
    \subseteq
    \mathcal{P}+\{\lfloor y_m^k\rfloor\mid k\ge1\},
    $$
    Romanoff's theorem shows that the latter set has positive lower asymptotic density for every $m\ge2$. Thus one obtains an explicit sequence of irrational algebraic bases tending to $3/2$ for which the qualitative conclusion is known, although the limiting base $3/2$ itself remains outside the reach of the present argument. This is consistent with the fixed base difficulty described in Tao's comment \cite{Tao244}.

    The conclusion of Theorem~\ref{thm-romanoff} does not imply density one. We give a concrete example in the interval $(1,2)$. The proof is related to covering congruences for Lucas numbers. Wang \cite{WangLucas} proved that a positive proportion of positive integers are not the sum of a Lucas number and a prime. In another paper, Wang \cite{WangAP} constructed arithmetic progressions avoiding sums involving Fibonacci and Lucas numbers.

    \begin{theorem}\label{thm-golden-obstruction}
        Let
        $$
        \varphi=\frac{1+\sqrt5}{2}.
        $$
        Then
        $$
        0<
        \liminf_{x\to\infty}
        \frac{\left|\mathcal{S}_\varphi\cap[1,x]\right|}{x}
        \le
        \limsup_{x\to\infty}
        \frac{\left|\mathcal{S}_\varphi\cap[1,x]\right|}{x}
        \le \frac{1937}{1938}.
        $$
        Moreover,
        $$
        \#\Bigl\{n\le x\mid n\notin \mathcal{P}+\bigl\{\lfloor \varphi^k\rfloor\mid k\ge1\bigr\}\Bigr\}
        \ge \frac{x}{1938}-O(\log x).
        $$
    \end{theorem}

    The obstruction in Theorem~\ref{thm-golden-obstruction} raises the question whether a base for which the complement has positive lower density can be found arbitrarily close to $1$. We conjecture that this is the case.

    \begin{conjecture}\label{conj-near-one-obstruction}
        For every $\varepsilon>0$, there exists a real number $a$ with $1<a<1+\varepsilon$ such that
        $$
        \liminf_{x\to\infty}\frac{1}{x}
        \#\Bigl\{n\le x\mid
        n\notin \mathcal{P}+\bigl\{\lfloor a^m\rfloor\mid m\ge1\bigr\}\Bigr\}>0.
        $$
    \end{conjecture}

    We prove a weaker assertion with infinitely many exceptional integers. For $a>1$, let
    \[
     E_a(x)=\#\bigl\{n\in\mathbb{N}\mid n\le x,\ n\notin\mathcal{S}_a\bigr\}.
    \]
    The following result holds for almost every base and gives a quantitative lower bound.

    \begin{theorem}\label{thm-prime-exceptions}
    For almost every real number $a>1$ and every $\eta>0$, there exists $x_0=x_0(a,\eta)$ such that
    \begin{equation}\label{eq-prime-exceptions}
     E_a(x)\ge x^{1-\eta}\qquad(x\ge x_0).
    \end{equation}
    In particular, for every $\varepsilon>0$, almost every $a\in(1,1+\varepsilon)$ has infinitely many positive integers that cannot be written as $p+\lfloor a^m\rfloor$ with $p$ prime and $m\ge1$.
    \end{theorem}

    Together with the trivial estimate $E_a(x)\le x$, the bounds in \eqref{eq-prime-exceptions} for all $\eta>0$ are equivalent to
    \[
     \lim_{x\to\infty}\frac{\log E_a(x)}{\log x}=1.
    \]
    It does not imply positive lower density, so Conjecture~\ref{conj-near-one-obstruction} remains stronger. Nevertheless, Theorems~\ref{thm-romanoff} and~\ref{thm-prime-exceptions} hold simultaneously for almost every $a>1$. Thus the same sumset has positive lower asymptotic density while its complement has logarithmic growth exponent one.

    The proof of Theorem~\ref{thm-prime-exceptions} is different from the Romanoff second moment argument. For almost every fixed $a$ and every sufficiently large integer $K$ that is a power of two, we construct a squarefree modulus $Q_K$ and an integer $b_K$ such that
    \[
     \gcd\bigl(b_K-\lfloor a^m\rfloor,Q_K\bigr)>1\qquad(1\le m\le K),
     \qquad
     \log Q_K=o(K).
    \]
    The congruences may change with $K$, but the base $a$ is fixed throughout. The subexponential size of $Q_K$ yields the lower bound in \eqref{eq-prime-exceptions}.

    Theorem~\ref{thm-prime-exceptions} leaves open the corresponding question for an arbitrary fixed base. It is natural to ask whether, for every real number $a>1$, there are infinitely many positive integers outside $\mathcal{S}_a$. This question is related to additive complements of the primes. Put
    \[
     B_a=\{\lfloor a^m\rfloor\mid m\ge1\}.
    \]
    Since the elements of $B_a$ are eventually distinct,
    \[
     \#(B_a\cap[1,x])=\frac{\log x}{\log a}+O_a(1).
    \]
    If, for some $a>1$, only finitely many positive integers lie outside $\mathcal{S}_a$, then $B_a$ is an additive complement to the primes with counting function $O(\log x)$. Ruzsa \cite[Theorem~2]{Ruzsa1998Primes} studied how thin such complements can be and proved that every additive complement $B$ to the primes satisfies
    \[
     \liminf_{x\to\infty}\frac{\#(B\cap[1,x])}{\log x}\ge e^\gamma,
    \]
    where $\gamma$ is Euler's constant. Consequently, every fixed $a$ with $\log a>e^{-\gamma}$ has infinitely many positive integers outside $\mathcal{S}_a$. This argument leaves the range $1<a\le \exp(e^{-\gamma})$ open.

    We next turn to squarefree integers. Erd\H{o}s \cite{Erdos1950} asked whether every sufficiently large odd integer is a squarefree integer plus a power of two. He returned to the problem in later lists \cite[p. 9]{Erdos1981}, and it also appears in Erd\H{o}s and Graham \cite[p. 28]{ErdosGraham}, in \cite{Erdos1997}, and as Problem 11 in Bloom's collection \cite{Bloom11}. There is a slight ambiguity in the historical record concerning the almost all version. Erd\H{o}s \cite{ErdosUnsolved1997} stated that he could prove that almost all integers not divisible by $4$ are the sum of a squarefree integer and a power of two, although we are not aware of a published proof. Granville and Soundararajan \cite{GranvilleSound} later showed that the original problem is closely related to the distribution of Wieferich primes. They proved that a positive answer to the original problem would imply that a positive proportion of primes $p$ satisfy $2^{p-1}\not\equiv1\pmod{p^2}$. In the other direction, under a quantitative sparsity hypothesis for the primes satisfying $2^{p-1}\equiv1\pmod{p^2}$, they proved that all but $O(x/\log x)$ odd integers not exceeding $x$ are the sum of a squarefree integer and a power of two. They also observed that the same conclusion follows from the sufficient condition
    $$
    \sum_{p^2\mid 2^{p-1}-1}\frac{1}{\operatorname{ord}_p(2)}\le\frac58.
    $$

    In the same paper, Granville and Soundararajan discussed the related problem concerning primes and powers of two. They recalled Erd\H{o}s's prediction that almost all odd integers should be the sum of a prime and two powers of two. Gallagher \cite{Gallagher} had proved that, for every $\varepsilon>0$, a suitable fixed number of powers of two gives a set of such sums with lower density at least $1-\varepsilon$. Motivated by this, Granville and Soundararajan asked whether one can prove unconditionally that almost all odd integers are the sum of a squarefree integer and at most $K$ powers of two for some fixed integer $K$. They proved this under the weaker summability assumption
    $$
    \sum_{p^2\mid 2^{p-1}-1}\frac{1}{\operatorname{ord}_p(2)}<\infty,
    $$
    and observed unconditionally that, for sufficiently large fixed $K$, all but $O(x/(K2^K))$ odd integers not exceeding $x$ have such a representation. Elsholtz, Planitzer and Schlage-Puchta have informed the author that they have obtained a positive answer to this unconditional almost all question. Hercher \cite{Hercher} verified the original conjecture with one power of two for every odd integer below $2^{50}$.

    The squarefree integers have asymptotic density $6/\pi^2$. They also form an asymptotic basis of order two, as follows from classical results discussed by Erd\H{o}s and Nathanson \cite[Introduction and Lemma 2]{ErdosNathanson}. The question above asks for a much thinner additive complement. Even the almost all version with a single power of two has no known published unconditional proof. Together with the weaker unconditional results available when several powers of two are allowed, this motivates the metric variant in which $2^m$ is replaced by $\lfloor a^m\rfloor$ for almost every real $a>1$. In this setting, the residue distribution of the sparse sequence becomes accessible to metric exponential sum estimates. Allowing two integer parts with distinct exponents gives a stronger exceptional set estimate. The following result gives both bounds.

    \begin{theorem}\label{thm-squarefree}
        Let $\mathcal{Q}$ denote the set of positive squarefree integers. For almost every real number $a>1$ and every $\varepsilon>0$,
        \begin{equation}\label{eq-squarefree-quantitative}
            \#\Bigl\{n\le x\mid n\notin \mathcal{Q}+\bigl\{\lfloor a^m\rfloor\mid m\ge1\bigr\}\Bigr\}
            \ll_{a,\varepsilon}
            \frac{x(\log\log x)^{1+\varepsilon}}{\sqrt{\log x}}.
        \end{equation}
        Moreover,
        \begin{equation}\label{eq-squarefree-two-powers}
            \#\Bigl\{n\le x\mid n\notin \mathcal{Q}+\bigl\{\lfloor a^m\rfloor+\lfloor a^\ell\rfloor\mid 1\le m<\ell\bigr\}\Bigr\}
            \ll_{a,\varepsilon}
            \frac{x(\log\log x)^{1+\varepsilon}}{\log x}.
        \end{equation}
        In particular, almost all positive integers can be represented in each of these two forms.
    \end{theorem}

    Theorem~\ref{thm-squarefree} is also related to additive complements of lacunary sequences. We refer to Ruzsa \cite{Ruzsa2001,Ruzsa2017} and to Fang and S\'andor \cite{FangSandor} for results in this direction. Relations between uniform distribution and the density of sumsets were studied by Volkmann \cite{Volkmann}. For recent results on discrete sumsets with one large summand, see Griesmer \cite{Griesmer}. The estimates in Theorem~\ref{thm-squarefree} follow from residue class estimates in which the square moduli are allowed to grow with the number of terms.

    The two problems above have different origins, but their proofs meet at the residue distribution of $\lfloor y^k\rfloor$. Koksma \cite{Koksma} proved that $y^k$ is uniformly distributed modulo one for almost all $y>1$. An exponential sum estimate gives, for almost all $y>1$,
    $$
    \#\{1\le k\le K\mid \lfloor y^k\rfloor\equiv r\pmod q\}
    =\frac Kq+o_{y,q}(K)
    $$
    for every $q\ge1$ and every $r\pmod q$. We also prove the exact weighted correlation formulas
    $$
    \sum_{\substack{1\le k<\ell\le K\\
            \lfloor y^k\rfloor\ne\lfloor y^\ell\rfloor}}
    W\!\left(\lfloor y^\ell\rfloor-\lfloor y^k\rfloor\right)
    =\left(\frac{\zeta(2)}{2\zeta(4)}+o_y(1)\right)K^2
    $$
    and
    $$
    \sum_{\substack{1\le k<\ell\le K\\
            \lfloor y^k\rfloor\ne\lfloor y^\ell\rfloor\\
            \lfloor y^k\rfloor\equiv\lfloor y^\ell\rfloor\pmod2}}
    W\!\left(\lfloor y^\ell\rfloor-\lfloor y^k\rfloor\right)
    =\left(\frac{3\zeta(2)}{10\zeta(4)}+o_y(1)\right)K^2,
    $$
    where
    $$
    W(n)=\prod_{p\mid n}\left(1+\frac1p\right).
    $$
    The second formula gives the sharper estimate used in the Romanoff second moment argument for Theorem~\ref{thm-romanoff}. The exponential sum estimate is used again in Section~\ref{sec-squarefree-proof} with growing square moduli to prove \eqref{eq-squarefree-quantitative}. The proof of \eqref{eq-squarefree-two-powers} uses an estimate with two parameters over separated blocks of exponents. For Theorem~\ref{thm-golden-obstruction}, the positive lower density follows from the theorem of Ballot and Luca applied to the odd indexed Lucas numbers, while the upper density bound follows from a finite covering congruence.

    For Theorem~\ref{thm-prime-exceptions}, we first prove joint residue estimates for finitely many separated exponents. We use these estimates to select groups of exponents whose integer parts have a common residue modulo a prime. A consequence of the theorem of Pippenger and Spencer \cite{PippengerSpencer}, stated for families of finite sets, makes the selections disjoint. Fourth moment estimates suffice because the construction is needed only when the number of exponents is a power of two. After covering the few remaining exponents by additional primes, the Chinese remainder theorem gives the required congruence class. The fact that the sum of the logarithms of the selected primes is $o(K)$ gives the estimate $\log Q_K=o(K)$ used above.

    Section~\ref{sec-romanoff-proof} proves the distribution estimates and Theorem~\ref{thm-romanoff}. Section~\ref{sec-golden-proof} proves Theorem~\ref{thm-golden-obstruction}. Section~\ref{sec-prime-exceptions} proves Theorem~\ref{thm-prime-exceptions}. Section~\ref{sec-squarefree-proof} proves the two quantitative squarefree estimates.

    \section{Preliminaries}

    Throughout the paper, $\mathbb{N}=\{1,2,3,\ldots\}$, $\mathcal{P}$ denotes the set of primes, and $\lfloor x\rfloor$ is the greatest integer not exceeding $x$. All assertions that hold for almost every point are with respect to Lebesgue measure, and we write $\operatorname{meas}$ for Lebesgue measure. We write $e(t)=\exp(2\pi i t)$. For a positive integer $n$, $\operatorname{rad}(n)$ denotes the product of the distinct prime divisors of $n$. For a set $\mathcal{A}\subset\mathbb{N}$, put
    $$
    \mathcal{A}(x)=\#(\mathcal{A}\cap[1,x]).
    $$
    Implicit constants may depend on a fixed compact interval of bases unless stated otherwise.

    We shall use the following form of the first lemma of Borel and Cantelli. See Rudin \cite[Theorem 1.41]{Rudin}.

    \begin{lemma}\label{lem-borel-cantelli}
        Let $\mathcal{E}_1,\mathcal{E}_2,\ldots$ be measurable sets. If
        $$
        \sum_{m\ge1}\operatorname{meas}(\mathcal{E}_m)<\infty,
        $$
        then almost every point belongs to only finitely many of the sets $\mathcal{E}_m$.
    \end{lemma}

    We shall also use the following form of the inequality of Erd\H{o}s and Tur\'an. See \cite[Ch. 2, Theorem 2.5]{KuipersNiederreiter}.

    \begin{lemma}\label{lem-erdos-turan}
        Let $x_1,\ldots,x_K$ be real numbers and let $I\subset[0,1)$ be an interval. For every $H\ge1$,
        $$
        \left|\#\bigl\{1\le k\le K\mid \{x_k\}\in I\bigr\}-K|I|\right|
        \le C\left(\frac{K}{H}+\sum_{1\le h\le H}\frac1h\left|\sum_{k\le K}e(hx_k)\right|\right),
        $$
        where $C$ is an absolute constant.
    \end{lemma}

    For a nonzero integer $h$, put
    $$
    W(h)=\prod_{p\mid h}\left(1+\frac1p\right).
    $$
    Define
    \begin{equation}\label{eq-sieve-constant}
        C_0=
        \limsup_{x\to\infty}
        \sup_{h\in\mathbb{Z}\setminus\{0\}}
        \frac{(\log x)^2}{xW(h)}
        \#\{p\le x\mid p\in\mathcal{P},\ p+h\in\mathcal{P}\}.
    \end{equation}
    The Selberg upper bound sieve in dimension two shows that the constant $C_0$ in \eqref{eq-sieve-constant} is finite. See Halberstam and Richert \cite[Ch. 5]{HalberstamRichert}. The definition of $C_0$ then gives the following uniform estimate.

    \begin{lemma}\label{lem-sieve}
        For every $\varepsilon>0$, uniformly in nonzero $h$ and all sufficiently large $x$,
        $$
        \#\{p\le x\mid p\in\mathcal{P},\ p+h\in\mathcal{P}\}
        \le (C_0+\varepsilon)\frac{x}{(\log x)^2}W(h).
        $$
    \end{lemma}

    \section{A Romanoff theorem for almost all real bases}\label{sec-romanoff-proof}

    \subsection{A second moment lemma}

    The following form of Romanoff's second moment argument will be used in the proof of Theorem~\ref{thm-romanoff}.

    \begin{lemma}\label{lem-romanoff-criterion}
        Let $\mathcal{A}=\{a_1<a_2<\cdots\}$ be an infinite subset of $\mathbb{N}$. Suppose that, for some $L>0$,
        \begin{equation}\label{eq-A-count-L}
            \mathcal{A}(x)=\frac{\log x}{L}+o(\log x)
        \end{equation}
        and
        \begin{equation}\label{eq-A-sum-small}
            \sum_{\substack{a\in\mathcal{A}\\a\le x}}a=O_{\mathcal{A}}(x).
        \end{equation}
        Assume also that
        \begin{equation}\label{eq-even-pair-condition}
            \sum_{\substack{a_1<a_2\le x\\
                    a_1,a_2\in\mathcal{A}\\
                    a_1\equiv a_2\pmod2}}
            W(a_2-a_1)
            \le \left(\frac{3\zeta(2)}{10\zeta(4)}+o(1)\right)\mathcal{A}(x)^2.
        \end{equation}
        Then
        $$
        \liminf_{x\to\infty}
        \frac{\#\{n\le x\mid n=p+a\text{ for some }p\in\mathcal{P},\ a\in\mathcal{A}\}}{x}
        \ge \frac{1}{L+9C_0/\pi^2}.
        $$
    \end{lemma}

    \begin{proof}
        For $n\le x$, let
        $$
        r_x(n)=\#\{(p,a)\mid p\in\mathcal{P},\ a\in\mathcal{A},\ a\le x/2,\ n=p+a\}.
        $$
        The prime number theorem, uniformly for $a\le x/2$, gives
        $$
        \pi(x-a)=\frac{x-a}{\log x}+o\left(\frac{x}{\log x}\right).
        $$
        It follows from \eqref{eq-A-count-L} and \eqref{eq-A-sum-small} that
        \begin{align}\label{eq-first-moment-sharp}
            \sum_{n\le x}r_x(n)
            =\sum_{\substack{a\in\mathcal{A}\\a\le x/2}}\pi(x-a)
            =\frac{1}{\log x}
            \sum_{\substack{a\in\mathcal{A}\\a\le x/2}}(x-a)+o(x)
            =\frac{x}{L}+o(x).
        \end{align}

        We next estimate the second moment. The diagonal contribution equals the first moment. Suppose that
        $$
        p_1+a_1=p_2+a_2,\qquad a_1<a_2.
        $$
        Then
        $$
        p_1-p_2=a_2-a_1.
        $$
        If $a_2-a_1$ is odd, one of $p_1,p_2$ must equal $2$. Since $p_1>p_2$, there is at most one such prime pair for each $(a_1,a_2)$. The total contribution of odd differences to the second moment is therefore
        $$
        O(\mathcal{A}(x)^2)=o(x).
        $$
        For even differences, Lemma~\ref{lem-sieve} shows that the two ordered pairs corresponding to $a_1$ and $a_2$ contribute at most
        $$
        2(C_0+\varepsilon)\frac{x}{(\log x)^2}
        W(a_2-a_1).
        $$
        Thus \eqref{eq-A-count-L} and \eqref{eq-even-pair-condition} give
        \begin{equation}\label{eq-second-moment-sharp}
            \sum_{n\le x}r_x(n)^2
            \le \frac{x}{L}
            +\frac35(C_0+\varepsilon)\frac{\zeta(2)}{\zeta(4)}\frac{x}{L^2}
            +o_{\varepsilon}(x).
        \end{equation}
        Since $\zeta(2)/\zeta(4)=15/\pi^2$, Cauchy's inequality, \eqref{eq-first-moment-sharp}, and \eqref{eq-second-moment-sharp} yield
        $$
        \liminf_{x\to\infty}
        \frac{\#\{n\le x\mid r_x(n)>0\}}x
        \ge \frac{1}{L+9(C_0+\varepsilon)/\pi^2}.
        $$
        Letting $\varepsilon$ tend to zero proves the lemma.
    \end{proof}

    \subsection{Metric estimates for residue classes}

    Fix $J=[Y_1,Y_2]$ with $1<Y_1<Y_2$. The next two lemmas give the exponential sum estimate and congruence collision bound required for Proposition~\ref{prop-metric}.

    \begin{lemma}\label{lem-exp-sum}
        Let $J=[Y_1,Y_2]$, where $1<Y_1<Y_2$, and let $B>0$ be fixed. If $K\ge2$, $1\le d\le K^B$, and $1\le h\le K$, then
        $$
        \int_J\left|\sum_{k\le K}e\left(\frac{h y^k}{d}\right)\right|^2\,\mathrm{d}y
        \ll_{Y_1,Y_2,B} K.
        $$
    \end{lemma}

    \begin{proof}
        Put $\alpha=h/d$. For $1\le k<\ell\le K$, set
        $$
        H_{k,\ell}(y)=y^\ell-y^k.
        $$
        For $y\in J$,
        $$
        H'_{k,\ell}(y)=\ell y^{\ell-1}-k y^{k-1}
        \ge (\ell-k)y^{\ell-1}
        \ge (\ell-k)Y_1^{\ell-1}.
        $$
        Moreover, $H'_{k,\ell}$ is increasing on $J$, since
        $$
        H''_{k,\ell}(y)=y^{k-2}\{\ell(\ell-1)y^{\ell-k}-k(k-1)\}>0
        $$
        for $y>1$ and $\ell>k$. The first derivative test therefore yields
        $$
        \left|\int_J e(\alpha H_{k,\ell}(y))\,\mathrm{d}y\right|
        \ll_{Y_1,Y_2}
        \min\left(1,\frac{d}{h(\ell-k)Y_1^{\ell-1}}\right).
        $$
        Expanding the square, it remains to estimate
        $$
        \sum_{1\le k<\ell\le K}
        \min\left(1,\frac{d}{h(\ell-k)Y_1^{\ell-1}}\right).
        $$
        If $d/h<1$, this sum is $O_{Y_1}(1)$. If $d/h\ge1$, put
        $$
        L_0=1+\Bigl\lfloor\frac{\log(d/h)}{\log Y_1}\Bigr\rfloor.
        $$
        The part with $\ell\le L_0+2$ is $O(L_0^2)$. For $\ell>L_0+2$ the remaining part is
        $$
        \ll \frac{d}{h}
        \sum_{\ell>L_0+2}\frac{1}{Y_1^{\ell-1}}
        \sum_{r=1}^{\ell-1}\frac1r
        \ll_{Y_1} L_0+1.
        $$
        Thus the off diagonal contribution is $O_{Y_1}(1+\log^2(2d/h))$, which is $O_{Y_1,B}(K)$ under $d\le K^B$ and $h\le K$. This proves the lemma.
    \end{proof}

    \begin{lemma}\label{lem-decreasing-weight}
        Let $I\subset\mathbb{R}$ be a bounded interval, let $d\ge1$, and let $g$ be a nonnegative nonincreasing function on $I$. Then
        $$
        \int_{I\cap\{u\mid \operatorname{dist}(u,d\mathbb{Z})<1\}} g(u)\,\mathrm{d}u
        \ll \frac1d\int_I g(u)\,\mathrm{d}u+\sup_I g,
        $$
        with an absolute implied constant.
    \end{lemma}

    \begin{proof}
        Partition $I$ into its nonempty intersections with intervals $[rd,(r+1)d)$. Within each such interval, the set of points at distance less than $1$ from $d\mathbb{Z}$ has length at most $2$. Bound the contributions of the first two intersections by $4\sup_I g$. For every subsequent intersection, the preceding interval of length $d$ is fully contained in $I$. Since $g$ is nonincreasing, its supremum on the current intersection is at most $d^{-1}$ times its integral over that preceding interval. The preceding intervals used in this way are disjoint. Summing gives the required bound.
    \end{proof}

    \begin{lemma}\label{lem-collision}
        Let $J=[Y_1,Y_2]$, where $1<Y_1<Y_2$. For $1\le k<\ell$ and $d\ge1$, define
        $$
        \mathcal{E}_{k,\ell}(d)=\bigl\{y\in J\mid
        \lfloor y^\ell\rfloor\equiv \lfloor y^k\rfloor\pmod d\bigr\}.
        $$
        Then
        $$
        \operatorname{meas}(\mathcal{E}_{k,\ell}(d))
        \ll_{Y_1,Y_2}
        \frac1d+\frac1{(\ell-k)Y_1^{\ell-1}}.
        $$
    \end{lemma}

    \begin{proof}
        Let $H(y)=y^\ell-y^k$. As in the proof of Lemma~\ref{lem-exp-sum}, $H$ and $H'$ are increasing on $J$, and
        $$
        H'(y)\ge (\ell-k)Y_1^{\ell-1}.
        $$
        If $\lfloor y^\ell\rfloor\equiv\lfloor y^k\rfloor\pmod d$, then
        $$
        \operatorname{dist}(H(y),d\mathbb{Z})<1.
        $$
        Make the change of variables $u=H(y)$. The density
        $$
        g(u)=\frac{1}{H'(H^{-1}(u))}
        $$
        is decreasing. By Lemma~\ref{lem-decreasing-weight},
        \begin{align*}
            \operatorname{meas}(\mathcal{E}_{k,\ell}(d))
            &\le
            \int_{H(J)\cap\{u\mid \operatorname{dist}(u,d\mathbb{Z})<1\}}g(u)\,\mathrm{d}u    \\
            &\ll \frac1d\int_{H(J)}g(u)\,\mathrm{d}u+g(H(Y_1)) \\
            &\ll_{Y_1,Y_2} \frac1d+\frac1{(\ell-k)Y_1^{\ell-1}}.
        \end{align*}
        This proves the lemma.
    \end{proof}

    Koksma \cite{Koksma} proved that $y^k$ is uniformly distributed modulo one for almost all $y>1$. We shall need the following residue class version.

    \begin{lemma}\label{lem-residue-ud}
        For almost every $y>1$, for every integer $q\ge1$ and every $0\le r<q$,
        $$
        \#\{1\le k\le K\mid \lfloor y^k\rfloor\equiv r\pmod q\}
        =\frac Kq+o_{y,q}(K)
        $$
        as $K\to\infty$.
    \end{lemma}

    \begin{proof}
        Fix $J=[Y_1,Y_2]$ with $1<Y_1<Y_2$, positive integers $q$ and $h$, and a real number $\lambda>1$. Put
        $$
        L_m=\lfloor \lambda^m\rfloor
        $$
        and
        $$
        S_K(y)=\sum_{k\le K}e\left(\frac{hy^k}{q}\right).
        $$
        By Lemma~\ref{lem-exp-sum}, for all sufficiently large $m$,
        $$
        \int_J |S_{L_m}(y)|^2\,\mathrm{d}y\ll_{J,q,h}L_m.
        $$
        For every positive integer $s$, Chebyshev's inequality gives
        $$
        \operatorname{meas}\bigl\{y\in J\mid |S_{L_m}(y)|>L_m/s\bigr\}
        \ll_{J,q,h}\frac{s^2}{L_m}.
        $$
        Since $L_m$ grows geometrically, $\sum_m L_m^{-1}<\infty$, and therefore
        $$
        \sum_{m=1}^{\infty}
        \operatorname{meas}\bigl\{y\in J\mid |S_{L_m}(y)|>L_m/s\bigr\}<\infty.
        $$
        By Lemma~\ref{lem-borel-cantelli}, and then taking a countable intersection over $s$, we obtain
        $$
        S_{L_m}(y)=o_{y,q,h,\lambda}(L_m)
        $$
        for almost all $y\in J$.

        If $L_m\le K<L_{m+1}$, then
        $$
        |S_K(y)|\le |S_{L_m}(y)|+L_{m+1}-L_m.
        $$
        Hence
        $$
        \limsup_{K\to\infty}\frac{|S_K(y)|}{K}\le \lambda-1.
        $$
        Taking the countable intersection over $\lambda=1+1/s$, and then over the positive integers $q$ and $h$, it follows that, for almost all $y\in J$,
        \begin{equation}\label{eq-fixed-mod-exp}
            \sum_{k\le K}e\left(\frac{hy^k}{q}\right)=o_{y,q,h}(K)
        \end{equation}
        for every fixed $q$ and $h$.

        Apply Lemma~\ref{lem-erdos-turan} to $x_k=y^k/q$ and the interval
        $$
        \left[\frac rq,\frac{r+1}{q}\right).
        $$
        For fixed $H$, divide by $K$, use \eqref{eq-fixed-mod-exp}, and let $K$ tend to infinity. This gives
        $$
        \limsup_{K\to\infty}
        \left|
        \frac1K\#\left\{k\le K\mid
        \left\{\frac{y^k}{q}\right\}\in
        \left[\frac rq,\frac{r+1}{q}\right)
        \right\}-\frac1q
        \right|
        \le \frac CH.
        $$
        Letting $H$ tend to infinity proves the assertion on $J$, because
        $$
        \left\{\frac{y^k}{q}\right\}\in
        \left[\frac rq,\frac{r+1}{q}\right)
        $$
        is equivalent to $\lfloor y^k\rfloor\equiv r\pmod q$. Taking $J=[1+1/n,n]$ for $n=2,3,\ldots$ completes the proof.
    \end{proof}

    For $K,d\ge1$ and $y>1$, define
    $$
    M_{K,d}(y)=\max_{0\le r<d}
    \#\{1\le k\le K\mid \lfloor y^k\rfloor\equiv r\pmod d\}.
    $$

    \begin{proposition}\label{prop-metric}
        Let $J=[Y_1,Y_2]$, where $1<Y_1<Y_2$. For almost all $y\in J$,
        \begin{equation}\label{eq-pair-asymptotic}
            \sum_{\substack{1\le k<\ell\le K\\
                    \lfloor y^k\rfloor\ne\lfloor y^\ell\rfloor}}
            W\!\left(\lfloor y^\ell\rfloor-\lfloor y^k\rfloor\right)
            =\left(\frac{\zeta(2)}{2\zeta(4)}+o_y(1)\right)K^2
        \end{equation}
        as $K\to\infty$.
    \end{proposition}

    \begin{proof}
        Choose $B>2$ so large that
        \begin{equation}\label{eq-B-choice}
            B\frac{\log Y_1}{\log Y_2}>3.
        \end{equation}
        The implied constants below may depend on $J$ and $B$.

        We begin with squarefree moduli not exceeding a power of the number of terms. Fix $\lambda>1$ and put
        $$
        L_m=\lfloor \lambda^m\rfloor,\qquad D_m=L_m^B.
        $$
        For $d\le D_m$, the congruence $\lfloor y^k\rfloor\equiv a_0\pmod d$ is equivalent to
        $$
        \left\{\frac{y^k}{d}\right\}\in
        \left[\frac{a_0}{d},\frac{a_0+1}{d}\right)
        $$
        for the representative $a_0\in\{0,1,\ldots,d-1\}$. Applying Lemma~\ref{lem-erdos-turan} with $H=L_m$ gives
        $$
        M_{L_m,d}(y)
        \le \frac{L_m}{d}
        +C
        +C\sum_{1\le h\le L_m}\frac1h
        \left|\sum_{k\le L_m}e\left(\frac{h y^k}{d}\right)\right|.
        $$
        Consequently
        \begin{equation}\label{eq-small-moduli-sharp}
            \sum_{\substack{d\le D_m\\ d\ \mathrm{squarefree}}}\frac{M_{L_m,d}(y)}d
            \le L_m\sum_{\substack{d\le D_m\\ d\ \mathrm{squarefree}}}\frac1{d^2}+O(\log L_m)+T_m(y),
        \end{equation}
        where
        $$
        T_m(y)=\sum_{\substack{d\le D_m\\ d\ \mathrm{squarefree}}}\sum_{h\le L_m}\frac1{dh}
        \left|\sum_{k\le L_m}e\left(\frac{h y^k}{d}\right)\right|.
        $$
        By Cauchy's inequality and Lemma~\ref{lem-exp-sum},
        \begin{align*}
            \int_J T_m(y)\,\mathrm{d}y\!
            \le
            \sum_{\substack{d\le D_m\\ d\ \mathrm{squarefree}}}\sum_{h\le L_m}\frac{1}{dh}
            \left(\operatorname{meas}(J)\int_J\left|\sum_{k\le L_m}e\left(\frac{h y^k}{d}\right)\right|^2\,\mathrm{d}y\right)^{1/2} 
            \!\ll L_m^{1/2}(\log L_m)^2.
        \end{align*}
        Hence, for each fixed $\eta>0$,
        $$
        \sum_{m=1}^{\infty}
        \operatorname{meas}\{y\in J\mid T_m(y)>\eta L_m\}<\infty.
        $$
        Applying Lemma~\ref{lem-borel-cantelli} with $\eta=1/s$, and then taking the countable intersection over positive integers $s$, gives, for almost all $y\in J$,
        $$
        T_m(y)=o_y(L_m).
        $$
        Since
        $$
        \sum_{\substack{d\ge1\\ d\ \mathrm{squarefree}}}\frac1{d^2}
        =\prod_p\left(1+\frac1{p^2}\right)
        =\frac{\zeta(2)}{\zeta(4)},
        $$
        it follows from \eqref{eq-small-moduli-sharp} that
        $$
        \sum_{\substack{d\le L_m^B\\ d\ \mathrm{squarefree}}}\frac{M_{L_m,d}(y)}d
        \le (\zeta(2)/\zeta(4)+o_y(1))L_m
        $$
        for almost all $y\in J$.

        We next pass from $L_m$ to arbitrary $K$. If $L_m\le K<L_{m+1}$, monotonicity gives
        $$
        \begin{aligned}
            \sum_{\substack{d\le K^B\\ d\ \mathrm{squarefree}}}\frac{M_{K,d}(y)}d
            \le \sum_{\substack{d\le L_{m+1}^B\\ d\ \mathrm{squarefree}}}\frac{M_{L_{m+1},d}(y)}d
            \le (\zeta(2)/\zeta(4)+o_y(1))L_{m+1}.
        \end{aligned}
        $$
        Since $L_{m+1}\le(\lambda+o(1))K$, taking the countable intersection over $\lambda=1+1/s$, where $s=1,2,\ldots$, gives
        \begin{equation}\label{eq-small-ae-sharp}
            \sum_{\substack{d\le K^B\\ d\ \mathrm{squarefree}}}\frac{M_{K,d}(y)}d
            \le (\zeta(2)/\zeta(4)+o_y(1))K.
        \end{equation}

        For each fixed $d$, let $m_a$ be the number of indices $k\le K$ for which $\lfloor y^k\rfloor\equiv a\pmod d$. Then
        \begin{align*}
            \#\{1\le k<\ell\le K\mid \lfloor y^k\rfloor\equiv\lfloor y^\ell\rfloor\pmod d\}
            &=\sum_a\binom{m_a}{2}\\
            &\le \frac12\sum_a m_aM_{K,d}(y)\\
            &=\frac K2M_{K,d}(y).
        \end{align*}
        After summing with weight $1/d$ over squarefree $d\le K^B$, \eqref{eq-small-ae-sharp} gives
        \begin{equation}\label{eq-small-pair-sharp}
            \sum_{\substack{1\le k<\ell\le K\\ \lfloor y^k\rfloor\ne\lfloor y^\ell\rfloor}}
            \sum_{\substack{d\mid \lfloor y^\ell\rfloor-\lfloor y^k\rfloor\\d\le K^B\\ d\ \mathrm{squarefree}}}\frac1d
            \le \left(\frac{\zeta(2)}{2\zeta(4)}+o_y(1)\right)K^2.
        \end{equation}

        It remains to consider divisors $d>K^B$. For the same sequence $L_m$, define
        $$
        R_m(y)=\sum_{1\le k<\ell\le L_{m+1}}
        \sum_{\substack{L_m^B<d\le Y_2^\ell\\
                \lfloor y^k\rfloor\ne\lfloor y^\ell\rfloor\\
                \lfloor y^\ell\rfloor\equiv\lfloor y^k\rfloor\pmod d}}
        \frac1d.
        $$
        If the inner congruence holds and the two floors are different, then
        $$
        d\le \left|\lfloor y^\ell\rfloor-\lfloor y^k\rfloor\right|\le Y_2^\ell.
        $$
        By Lemma~\ref{lem-collision},
        \begin{align*}
            \int_J R_m(y)\,\mathrm{d}y
            &\ll \sum_{1\le k<\ell\le L_{m+1}}
            \sum_{L_m^B<d\le Y_2^\ell}
            \left(\frac1{d^2}+\frac1{d(\ell-k)Y_1^{\ell-1}}\right)\\
            &\ll \frac{L_{m+1}^2}{L_m^B}
            +\sum_{\substack{\ell\le L_{m+1}\\Y_2^\ell>L_m^B}}
            \frac{1+\log(Y_2^\ell/L_m^B)}{Y_1^{\ell-1}}
            \sum_{r=1}^{\ell-1}\frac1r.
        \end{align*}
        Put $\ell_0=B\log L_m/\log Y_2$. The second term is
        $$
        \ll_J \sum_{\ell>\ell_0}\frac{\ell\log \ell}{Y_1^\ell}
        \ll_J (\ell_0+1)^2Y_1^{-\ell_0}.
        $$
        By \eqref{eq-B-choice}, this is $O(L_m^{-3})$, while
        $$
        L_{m+1}^2/L_m^B\ll_\lambda L_m^{2-B}.
        $$
        Consequently, for every fixed $\eta>0$,
        $$
        \sum_{m=1}^{\infty}
        \operatorname{meas}\{y\in J\mid R_m(y)>\eta L_m^2\}<\infty.
        $$
        Applying Lemma~\ref{lem-borel-cantelli} with $\eta=1/s$, and then taking the countable intersection over positive integers $s$, gives
        $$
        R_m(y)=o_y(L_m^2)
        $$
        for almost all $y\in J$. If $L_m\le K<L_{m+1}$, every contribution with $k,\ell\le K$ and $d>K^B$ is included in $R_m(y)$. Hence
        \begin{equation}\label{eq-large-pair-negligible}
            \sum_{\substack{1\le k<\ell\le K\\ \lfloor y^k\rfloor\ne\lfloor y^\ell\rfloor}}
            \sum_{\substack{d\mid \lfloor y^\ell\rfloor-\lfloor y^k\rfloor\\d>K^B}}\frac1d
            =o_y(K^2).
        \end{equation}
        Finally,
        $$
        W(n)=\sum_{d\mid \operatorname{rad}(n)}\frac{1}{d}
        $$
        for every positive integer $n$. Combining \eqref{eq-small-pair-sharp} and \eqref{eq-large-pair-negligible} proves
        \begin{equation}\label{eq-pair-upper}
            \limsup_{K\to\infty}\frac1{K^2}
            \sum_{\substack{1\le k<\ell\le K\\
                    \lfloor y^k\rfloor\ne\lfloor y^\ell\rfloor}}
            W\!\left(\lfloor y^\ell\rfloor-\lfloor y^k\rfloor\right)
            \le \frac{\zeta(2)}{2\zeta(4)}.
        \end{equation}

        For the reverse inequality, intersect with the set of full measure on which Lemma~\ref{lem-residue-ud} holds. Fix $D\ge1$. For a squarefree integer $d\le D$ and $0\le r<d$, put
        $$
        N_r(K)=\#\{k\le K\mid \lfloor y^k\rfloor\equiv r\pmod d\}.
        $$
        By Lemma~\ref{lem-residue-ud}, for almost all $y\in J$,
        $$
        N_r(K)=\frac Kd+o_{y,d}(K).
        $$
        Hence, for every fixed squarefree $d$,
        $$
        C_d(K)=\sum_{r=0}^{d-1}\binom{N_r(K)}2
        $$
        satisfies
        $$
        C_d(K)=\frac12\sum_{r=0}^{d-1}N_r(K)^2-\frac K2
        $$
        and therefore
        $$
        C_d(K)=\frac{K^2}{2d}+o_{y,d}(K^2).
        $$
        Since $y^{k+1}-y^k=(y-1)y^k>1$ for all sufficiently large $k$, the sequence $\lfloor y^k\rfloor$ is eventually strictly increasing. Hence the number of pairs $k<\ell$ with $\lfloor y^k\rfloor=\lfloor y^\ell\rfloor$ is bounded in terms of $y$, and
        $$
        \#\{k<\ell\le K\mid \lfloor y^k\rfloor\ne\lfloor y^\ell\rfloor,\ d\mid \lfloor y^\ell\rfloor-\lfloor y^k\rfloor\}
        =\frac{K^2}{2d}+o_{y,d}(K^2).
        $$
        Keeping only the squarefree divisors not exceeding $D$ in the divisor expansion of $W$ gives
        $$
        \begin{aligned}
            &\sum_{\substack{k<\ell\le K\\ \lfloor y^k\rfloor\ne\lfloor y^\ell\rfloor}}
            W\!\left(\lfloor y^\ell\rfloor-\lfloor y^k\rfloor\right)\\
            &\qquad\qquad\ge
            \sum_{\substack{d\le D\\d\ \mathrm{squarefree}}}\frac1d
            \#\{k<\ell\le K\mid \lfloor y^k\rfloor\ne\lfloor y^\ell\rfloor,\ d\mid \lfloor y^\ell\rfloor-\lfloor y^k\rfloor\}.
        \end{aligned}
        $$
        Consequently
        $$
        \liminf_{K\to\infty}\frac1{K^2}
        \sum_{\substack{k<\ell\le K\\ \lfloor y^k\rfloor\ne\lfloor y^\ell\rfloor}}
        W\!\left(\lfloor y^\ell\rfloor-\lfloor y^k\rfloor\right)
        \ge
        \frac12\sum_{\substack{d\le D\\d\ \mathrm{squarefree}}}\frac1{d^2}.
        $$
        Letting $D$ tend to infinity gives
        $$
        \liminf_{K\to\infty}\frac1{K^2}
        \sum_{\substack{k<\ell\le K\\ \lfloor y^k\rfloor\ne\lfloor y^\ell\rfloor}}
        W\!\left(\lfloor y^\ell\rfloor-\lfloor y^k\rfloor\right)
        \ge \frac{\zeta(2)}{2\zeta(4)}.
        $$
        Together with \eqref{eq-pair-upper}, this gives \eqref{eq-pair-asymptotic}.
    \end{proof}

    We shall also need the corresponding asymptotic restricted to even differences.

    \begin{proposition}\label{prop-metric-even}
        For almost every $y>1$,
        \begin{equation}\label{eq-even-pair-asymptotic}
            \sum_{\substack{1\le k<\ell\le K\\
                    \lfloor y^k\rfloor\ne\lfloor y^\ell\rfloor\\
                    \lfloor y^k\rfloor\equiv\lfloor y^\ell\rfloor\pmod2}}
            W\!\left(\lfloor y^\ell\rfloor-\lfloor y^k\rfloor\right)
            =\left(\frac{3\zeta(2)}{10\zeta(4)}+o_y(1)\right)K^2.
        \end{equation}
    \end{proposition}

    \begin{proof}
        It is enough to work on a fixed interval $J=[Y_1,Y_2]$ with $1<Y_1<Y_2$. If $n$ is even, then
        \begin{equation}\label{eq-even-W-expansion}
            W(n)=\frac32
            \sum_{\substack{d\mid n\\d\ \mathrm{odd,\ squarefree}}}\frac1d.
        \end{equation}
        We first obtain the upper bound. Fix $B>2$ satisfying \eqref{eq-B-choice}. The proof of \eqref{eq-large-pair-negligible}, with this same $B$, applies unchanged, and we intersect with the corresponding set of full measure. Fix $\lambda>1$, put $L_m=\lfloor \lambda^m\rfloor$, and set $D_m=L_m^B$. Lemma~\ref{lem-erdos-turan}, applied modulo $2d$, gives
        $$
        M_{L_m,2d}(y)
        \le
        \frac{L_m}{2d}+C
        +C\sum_{h\le L_m}\frac1h
        \left|\sum_{k\le L_m}e\left(\frac{h y^k}{2d}\right)\right|.
        $$
        After multiplication by $1/d$ and summation over odd squarefree $d\le D_m$, this becomes
        $$
        \sum_{\substack{d\le D_m\\d\ \mathrm{odd,\ squarefree}}}
        \frac{M_{L_m,2d}(y)}d
        \le
        \frac{L_m}{2}
        \sum_{\substack{d\le D_m\\d\ \mathrm{odd,\ squarefree}}}\frac1{d^2}
        +O(\log L_m)+T_m^{(2)}(y),
        $$
        where
        $$
        T_m^{(2)}(y)=
        \sum_{\substack{d\le D_m\\d\ \mathrm{odd,\ squarefree}}}
        \sum_{h\le L_m}\frac1{dh}
        \left|\sum_{k\le L_m}e\left(\frac{h y^k}{2d}\right)\right|.
        $$
        Cauchy's inequality and Lemma~\ref{lem-exp-sum}, applied with exponent $B+1$, give
        $$
        \int_J T_m^{(2)}(y)\,\mathrm{d}y
        \ll L_m^{1/2}(\log L_m)^2.
        $$
        Lemma~\ref{lem-borel-cantelli} therefore gives $T_m^{(2)}(y)=o_y(L_m)$ for almost all $y\in J$. Since
        $$
        \sum_{\substack{d\ge1\\d\ \mathrm{odd,\ squarefree}}}\frac1{d^2}
        =
        \prod_{p>2}\left(1+\frac1{p^2}\right)
        =
        \frac45\frac{\zeta(2)}{\zeta(4)},
        $$
        we obtain along the sequence $L_m$
        $$
        \sum_{\substack{d\le L_m^B\\d\ \mathrm{odd,\ squarefree}}}
        \frac{M_{L_m,2d}(y)}d
        \le
        \left(\frac{2\zeta(2)}{5\zeta(4)}+o_y(1)\right)L_m.
        $$
        The monotonicity argument used in the proof of \eqref{eq-small-ae-sharp}, followed by a countable intersection over $\lambda\downarrow1$, gives
        \begin{equation}\label{eq-even-small-moduli}
            \sum_{\substack{d\le K^B\\d\ \mathrm{odd,\ squarefree}}}
            \frac{M_{K,2d}(y)}d
            \le
            \left(\frac{2\zeta(2)}{5\zeta(4)}+o_y(1)\right)K.
        \end{equation}

        For each odd squarefree $d$, the number of pairs $k<\ell\le K$ with
        $$
        2d\mid \lfloor y^\ell\rfloor-\lfloor y^k\rfloor
        $$
        is at most
        $$
        \frac K2 M_{K,2d}(y).
        $$
        It follows from \eqref{eq-even-W-expansion} and \eqref{eq-even-small-moduli} that the contribution of $d\le K^B$ is at most
        $$
        \left(\frac{3\zeta(2)}{10\zeta(4)}+o_y(1)\right)K^2.
        $$
        For the same choice of $B$, the contribution of $d>K^B$ is $o_y(K^2)$ by \eqref{eq-large-pair-negligible}. This proves the required upper bound.

        For the reverse inequality, fix an odd squarefree integer $d$. Lemma~\ref{lem-residue-ud}, applied modulo $2d$, gives
        $$
        \#\{k<\ell\le K\mid
        2d\mid \lfloor y^\ell\rfloor-\lfloor y^k\rfloor\}
        =\frac{K^2}{4d}+o_{y,d}(K^2).
        $$
        As in the proof of Proposition~\ref{prop-metric}, pairs with equal integer parts contribute only $O_y(1)$. Keeping in \eqref{eq-even-W-expansion} only odd squarefree divisors $d\le D$, and then letting $K$ tend to infinity, gives
        $$
        \liminf_{K\to\infty}\frac1{K^2}
        \sum_{\substack{k<\ell\le K\\
                \lfloor y^k\rfloor\ne\lfloor y^\ell\rfloor\\
                \lfloor y^k\rfloor\equiv\lfloor y^\ell\rfloor\pmod2}}
        W\!\left(\lfloor y^\ell\rfloor-\lfloor y^k\rfloor\right)
        \ge
        \frac38
        \sum_{\substack{d\le D\\d\ \mathrm{odd,\ squarefree}}}\frac1{d^2}.
        $$
        Finally,
        $$
        \sum_{\substack{d\ge1\\d\ \mathrm{odd,\ squarefree}}}\frac1{d^2}
        =\prod_{p>2}\left(1+\frac1{p^2}\right)
        =\frac45\frac{\zeta(2)}{\zeta(4)}.
        $$
        Letting $D$ tend to infinity proves \eqref{eq-even-pair-asymptotic}.
    \end{proof}

    \subsection{Proof of Theorem~\ref{thm-romanoff}}

    \begin{proof}[Proof of Theorem~\ref{thm-romanoff}]
        It is enough to work on a fixed compact interval $J=[Y_1,Y_2]\subset(1,\infty)$. Proposition~\ref{prop-metric-even} gives, for almost all $y\in J$,
        \begin{equation}\label{eq-even-pair-for-indices}
            \sum_{\substack{1\le k<\ell\le K\\
                    \lfloor y^k\rfloor\ne\lfloor y^\ell\rfloor\\
                    \lfloor y^k\rfloor\equiv\lfloor y^\ell\rfloor\pmod2}}
            W\!\left(\lfloor y^\ell\rfloor-\lfloor y^k\rfloor\right)
            \le
            \left(\frac{3\zeta(2)}{10\zeta(4)}+o_y(1)\right)K^2.
        \end{equation}
        Fix such a $y$.

        Let
        $$
        \mathcal{B}_y=\{\lfloor y^k\rfloor\mid k\in\mathbb{N}\},
        $$
        where repeated values are removed. Since $\lfloor y^k\rfloor$ is eventually strictly increasing,
        $$
        \mathcal{B}_y(x)=\frac{\log x}{\log y}+o_y(\log x).
        $$
        Also
        $$
        \sum_{\substack{b\in\mathcal{B}_y\\b\le x}}b=O_y(x),
        $$
        because the terms are bounded by a geometric progression up to an additive constant.

        It remains to verify \eqref{eq-even-pair-condition}. Let $K$ be the largest index for which $\lfloor y^K\rfloor\le x$. Every pair of distinct values $b_1<b_2\le x$ in $\mathcal{B}_y$ occurs as $b_1=\lfloor y^k\rfloor$ and $b_2=\lfloor y^\ell\rfloor$ for some $k<\ell\le K$. Hence the sum over distinct values is bounded above by the corresponding indexed sum in \eqref{eq-even-pair-for-indices}. Since $\lfloor y^k\rfloor$ is eventually strictly increasing, $K=\mathcal{B}_y(x)+O_y(1)$. Therefore \eqref{eq-even-pair-for-indices} yields
        $$
        \sum_{\substack{b_1<b_2\le x\\
                b_1,b_2\in\mathcal{B}_y\\
                b_1\equiv b_2\pmod2}}
        W(b_2-b_1)
        \le
        \left(\frac{3\zeta(2)}{10\zeta(4)}+o_y(1)\right)\mathcal{B}_y(x)^2.
        $$
        Lemma~\ref{lem-romanoff-criterion}, with $L=\log y$, now gives
        $$
        \delta_y\ge \frac{1}{\log y+9C_0/\pi^2}.
        $$
        This holds for almost all $y$ in every compact subinterval of $(1,\infty)$, and hence for almost all $y>1$.
    \end{proof}

    \section{The golden ratio}\label{sec-golden-proof}

    We prove Theorem~\ref{thm-golden-obstruction}. From this point on, $L_m$ denotes the Lucas numbers. Let
    $$
    L_0=2,\qquad L_1=1,\qquad L_{m+2}=L_{m+1}+L_m
    $$
    be the Lucas numbers. Binet's formula gives
    $$
    L_m=\varphi^m+(-\varphi^{-1})^m.
    $$
    Since $0<\varphi^{-m}<1$ for $m\ge1$, if
    $$
    u_m=\lfloor \varphi^m\rfloor,
    $$
    then
    \begin{equation}\label{eq-floor-phi-lucas}
        u_m=
        \begin{cases}
            L_m, & m\text{ odd},\\
            L_m-1, & m\text{ even}.
        \end{cases}
    \end{equation}
    The same formula is also valid for $m=0$.

    We first prove the positive lower density assertion. For $j\ge0$, put
    $$
    v_j=L_{2j+1}.
    $$
    By \eqref{eq-floor-phi-lucas},
    $$
    v_j=\lfloor\varphi^{2j+1}\rfloor.
    $$
    The sequence $(v_j)$ satisfies
    $$
    v_{j+2}=3v_{j+1}-v_j,
    $$
    and its characteristic polynomial is $X^2-3X+1$, with distinct roots $\varphi^2$ and $\varphi^{-2}$ whose quotient $\varphi^4$ is not a root of unity. It is therefore a nonconstant nondegenerate integral linear recurrence with separable characteristic polynomial. By the theorem of Ballot and Luca \cite[Theorem~1]{BallotLuca}, the set
    $$
    \mathcal{P}+\{v_j\mid j\ge0\}
    $$
    has positive lower asymptotic density. Since this set is contained in $\mathcal{S}_\varphi$, we obtain
    $$
    \liminf_{x\to\infty}
    \frac{\left|\mathcal{S}_\varphi\cap[1,x]\right|}{x}>0.
    $$

    It remains to prove the upper density assertion. The following table contains three coverings. The first three rows are used in all three. They are completed by the rows marked $\mathrm{I}$, $\mathrm{II}$, or $\mathrm{III}$. In each row, $u_m$ is periodic modulo $\ell$ with period dividing $T$, and the last column gives the residue classes of $m\pmod T$ for which $u_m\equiv c\pmod\ell$.
    $$
    \begin{array}{c|c|c|c|l}
        & \ell & c & T & m\pmod T \\
        \hline
        & 2  & 1  & 6  & 0,1,5\\
        & 3  & 2  & 8  & 2,5,6,7\\
        & 5  & 1  & 4  & 0,1\\
        \mathrm{I}  & 7  & 3  & 16 & 6,10,11,13\\
        \mathrm{I}  & 23 & 4  & 48 & 3,18,21,30\\
        \mathrm{II} & 7  & 4  & 16 & 3,5,8\\
        \mathrm{II} & 23 & 19 & 48 & 22,26,27,45\\
        \mathrm{III}& 17 & 4  & 36 & 3,15\\
        \mathrm{III}& 19 & 0  & 18 & 9
    \end{array}
    $$
    For every row,
    $$
    L_T\equiv L_0\pmod\ell,
    \qquad
    L_{T+1}\equiv L_1\pmod\ell.
    $$
    The recurrence therefore shows that $L_m$ has period dividing $T$ modulo $\ell$. Since every listed $T$ is even, \eqref{eq-floor-phi-lucas} gives the same period for $u_m$. Direct substitution verifies the classes in the last column. For the first two coverings, the first three rows cover every residue class modulo $48$ except $3$ and $27$. In covering $\mathrm{I}$, these two classes are supplied by the row for $23$ and the row for $7$, respectively. In covering $\mathrm{II}$, they are supplied by the row for $7$ and the row for $23$, respectively. For covering $\mathrm{III}$, work modulo $72$. The first three rows leave only the classes $3,27,51$. The row for $17$ covers $3$ and $51$, and the row for $19$ covers $27$. Thus each of the three choices covers all exponent classes.

    The three corresponding systems of congruences are summarized by
    $$
    \begin{array}{c|ccccccc}
        & 2 & 3 & 5 & 7 & 17 & 19 & 23\\
        \hline
        1361 & 1 & 2 & 1 & 3 &   &   & 4\\
        3791 & 1 & 2 & 1 & 4 &   &   & 19\\
        6821 & 1 & 2 & 1 &   & 4 & 0 &
    \end{array}
    $$
    where each displayed entry is the residue modulo the prime at the top of its column. The first two rows are read modulo
    $$
    2\cdot3\cdot5\cdot7\cdot23=4830,
    $$
    and the third is read modulo
    $$
    2\cdot3\cdot5\cdot17\cdot19=9690.
    $$
    Let
    $$
    \begin{aligned}
        \mathcal{N}={}&\{n\in\mathbb{N}\mid n\equiv1361\pmod{4830}
        \text{ or }n\equiv3791\pmod{4830}\}\\
        &{}\cup\{n\in\mathbb{N}\mid n\equiv6821\pmod{9690}\}.
    \end{aligned}
    $$
    Fix $n\in\mathcal{N}$ and $m\ge1$. Choose one of the three congruence classes above that contains $n$, and use the corresponding covering. There is a row $(\ell,c,T)$ for which
    $$
    u_m\equiv c\pmod\ell
    \qquad\text{and}\qquad
    n\equiv c\pmod\ell.
    $$
    Hence $\ell\mid n-u_m$. If $n=p+u_m$ with $p$ prime, then $p=\ell$. Therefore every represented element of $\mathcal{N}$ belongs to
    $$
    \{\ell+u_m\mid m\ge1,\ \ell\in\{2,3,5,7,17,19,23\}\}.
    $$
    There are only $O(\log x)$ such integers not exceeding $x$. Moreover,
    $$
    \gcd(4830,9690)=30
    $$
    and
    $$
    1361\equiv3791\equiv6821\equiv11\pmod {30}.
    $$
    Thus each of the first two progressions meets the third in one residue class modulo
    $$
    \operatorname{lcm}(4830,9690)=1560090.
    $$
    The two intersections are disjoint, and hence
    $$
    \begin{aligned}
        \#(\mathcal{N}\cap[1,x])
        &=\left(\frac{2}{4830}+\frac{1}{9690}-\frac{2}{1560090}\right)x+O(1)\\
        &=\frac{x}{1938}+O(1).
    \end{aligned}
    $$
    Consequently
    $$
    \#\Bigl\{n\le x\mid n\notin \mathcal{P}+\bigl\{\lfloor \varphi^k\rfloor\mid k\ge1\bigr\}\Bigr\}
    \ge \frac{x}{1938}-O(\log x).
    $$
    Hence
    $$
    \limsup_{x\to\infty}
    \frac{\left|\mathcal{S}_\varphi\cap[1,x]\right|}{x}
    \le\frac{1937}{1938},
    $$
    which completes the proof of Theorem~\ref{thm-golden-obstruction}.

    \section{Integers not represented by a prime and a real power}\label{sec-prime-exceptions}

We prove Theorem~\ref{thm-prime-exceptions}. The congruences used in this section may depend on the number of exponents under consideration. This is different from the fixed congruences in Section~\ref{sec-golden-proof}.

For $a>1$ and $m\ge1$, write $u_m(a)=\lfloor a^m\rfloor$. We shall prove the following proposition. It is enough to construct the congruences when the number of exponents is a power of two.

\begin{proposition}\label{prop-prime-cover}
For almost every $a>1$ and every sufficiently large integer $K$ that is a power of two, there are a squarefree integer $Q_K\ge2$ and an integer $b_K$ such that
\begin{equation}\label{eq-prime-cover}
 \gcd\bigl(b_K-u_m(a),Q_K\bigr)>1
 \qquad(1\le m\le K),
\end{equation}
all prime divisors of $Q_K$ are at most $8K$, and
\begin{equation}\label{eq-prime-cover-size}
 \log Q_K=o(K)
\end{equation}
as $K$ tends to infinity through powers of two.
\end{proposition}

We first explain why this proposition gives the required exceptional integers.

\begin{proof}[Proof of Theorem~\ref{thm-prime-exceptions} from Proposition~\ref{prop-prime-cover}]
Fix a base $a$ satisfying the proposition. For a sufficiently large real number $x$, put
\[
 M=\left\lceil\frac{\log(x+1)}{\log a}\right\rceil
\]
and let $K$ be the smallest power of two with $K\ge M$. Then
\[
 M\le K<2M.
\]
If $m>M$, then $a^m>x+1$, and hence $u_m(a)>x$. Such an exponent cannot occur in a representation of a positive integer not exceeding $x$.

Apply Proposition~\ref{prop-prime-cover} at the dyadic integer $K$. For each $m\le M$, choose a prime $q_m$ dividing both $Q_K$ and $b_K-u_m(a)$, which is possible by \eqref{eq-prime-cover}. Suppose that $n\equiv b_K\pmod{Q_K}$ and $n=p+u_m(a)$ for a prime $p$. Then
\[
 p=n-u_m(a)\equiv b_K-u_m(a)\equiv0\pmod{q_m}.
\]
Since $p$ is prime, this implies $p=q_m$. Thus every represented integer in this progression belongs to
\[
 \{u_m(a)+q_m\mid1\le m\le M\},
\]
which contains at most $M\le K$ integers. There are at least $x/Q_K-1$ positive integers not exceeding $x$ in any residue class modulo $Q_K$. Therefore
\begin{equation}\label{eq-prime-progression-count}
 E_a(x)\ge \frac{x}{Q_K}-K-1.
\end{equation}

We have $M=(\log x)/(\log a)+O_a(1)$ and $K<2M$. Hence \eqref{eq-prime-cover-size} gives $\log Q_K=o(\log x)$. Fix $0<\eta<1$. For all sufficiently large $x$,
\[
 Q_K\le x^{\eta/2},\qquad
 K+1\le \tfrac12x^{1-\eta/2},\qquad
 \tfrac12x^{\eta/2}\ge1.
\]
Substitution in \eqref{eq-prime-progression-count} gives
\[
 E_a(x)\ge \tfrac12x^{1-\eta/2}\ge x^{1-\eta}.
\]
The assertion for $\eta\ge1$ follows from the case $\eta=1/2$. This proves \eqref{eq-prime-exceptions}. In particular, the set of bases satisfying this estimate has full measure in every interval $(1,1+\varepsilon)$, and the estimate with $\eta=1/2$ gives infinitely many exceptions for each such base.
\end{proof}

\subsection{Joint distribution of separated powers}

Fix a compact interval $I=[A,B]$ with $1<A<B$. In this subsection, $a$ is chosen uniformly from $I$. Thus, for a measurable set $F\subseteq I$, we write
\[
 \mathbb P_I(F)=\frac{\operatorname{meas}(F)}{B-A},
 \qquad
 \mathbb E_I f=\frac1{B-A}\int_A^B f(a)\,\mathrm{d}a.
\]
The first quantity is the probability of $F$, and the second is the expectation of $f(a)$. An event involving $a$ is identified with the set of bases for which it holds.

We need a joint estimate for finitely many exponents that are sufficiently far apart. For an integer $K\ge3$, put
\[
 d=\lceil(\log K)^4\rceil.
\]
We shall only use values of $K$ large enough that $d<K$.

\begin{lemma}\label{lem-prime-joint-residues}
Fix a positive integer $t$ and a real number $C>0$. Suppose that
\begin{equation}\label{eq-prime-separated-indices}
 d\le m_1<\cdots<m_t\le K,
 \qquad m_{i+1}-m_i\ge d\quad(1\le i<t),
\end{equation}
and that $1\le M_i\le\exp(C(\log K)^2)$ are integers. For every set
\[
 F\subseteq\prod_{i=1}^t\{0,1,\ldots,M_i-1\},
\]
we have
\begin{equation}\label{eq-prime-joint-residue-error}
 \left|
 \mathbb P_I\bigl((u_{m_i}(a)\bmod M_i)_{i=1}^t\in F\bigr)
 -\frac{|F|}{M_1\cdots M_t}
 \right|
 \le \exp(-c_I d)
\end{equation}
for all sufficiently large $K$. The constant $c_I>0$ depends only on $I$. The threshold for $K$ may depend on $I,t,C$, but not on the exponents, the moduli, or $F$.
\end{lemma}

\begin{proof}
For an integer $M\ge1$ and a residue $r\pmod M$, define
\[
 h_{M,r}(y)=\mathbf1_{\{\lfloor y\rfloor\equiv r\pmod M\}}-\frac1M,
\]
where an indicator is equal to $1$ when its stated condition holds and to $0$ otherwise. The function $h_{M,r}$ has period $M$ and integral zero over one period. On each period its positive part is supported on an interval of length $1$. It follows that
\[
 H_{M,r}(y)=\int_0^y h_{M,r}(v)\,\mathrm{d}v
\]
is periodic and satisfies $|H_{M,r}(y)|\le1$ for all real $y$.

Let $A\le \alpha<\beta\le B$ and $n\ge1$. Integration by parts gives
\[
 \int_\alpha^\beta h_{M,r}(a^n)\,\mathrm{d}a
 =\left[\frac{H_{M,r}(a^n)}{na^{n-1}}\right]_\alpha^\beta
  -\int_\alpha^\beta H_{M,r}(a^n)
       \left(\frac1{na^{n-1}}\right)'\,\mathrm{d}a.
\]
The absolute value of the boundary term is at most $2/(n\alpha^{n-1})$. Since $1/(na^{n-1})$ is nonincreasing, the absolute value of the integral on the right is at most $1/(n\alpha^{n-1})$. Hence
\begin{equation}\label{eq-prime-periodic-integral}
 \left|\int_\alpha^\beta h_{M,r}(a^n)\,\mathrm{d}a\right|
 \le\frac3{n\alpha^{n-1}}.
\end{equation}
The integrands have only finitely many discontinuities on a compact interval, so this integration by parts can also be performed separately between consecutive discontinuities.

Fix residues $r_i\pmod{M_i}$. We first estimate the probability of this one residue choice. Suppose $t\ge2$, let $n=m_t$, and put
\[
 g(a)=\prod_{i=1}^{t-1}
       \mathbf1_{\{u_{m_i}(a)\equiv r_i\pmod{M_i}\}}.
\]
Partition $I$ into intervals $J=[\alpha,\beta]$ such that
\[
 \log(\beta/\alpha)\le\frac{d\log A}{4K}.
\]
This can be done with $O_I(K/d+1)$ intervals by dividing $[\log A,\log B]$ into equal intervals of the indicated maximum length.

The value of $g$ can change only when $a^{m_i}$ crosses an integer for some $i<t$. On $J$, the number of such crossings is at most
\[
 \sum_{i<t}(\beta^{m_i}+2)
 \le (t-1)(\beta^{n-d}+2),
\]
where we used \eqref{eq-prime-separated-indices}. On every interval between consecutive crossings, $g$ is constant and belongs to $\{0,1\}$. Applying \eqref{eq-prime-periodic-integral} on the intervals where $g=1$ gives
\[
 \left|\int_J g(a)h_{M_t,r_t}(a^n)\,\mathrm{d}a\right|
 \ll_t\frac{\beta^{n-d}+1}{n\alpha^{n-1}}.
\]
Furthermore,
\[
 \frac{\beta^{n-d}}{\alpha^{n-1}}
 =\alpha^{1-d}(\beta/\alpha)^{n-d}
 \le B A^{-d} A^{d/4}=B A^{-3d/4},
\]
and $\alpha^{1-n}\le B A^{-d}$. After summing over the partition of $I$, we obtain
\[
 \int_I g(a)\mathbf1_{\{u_n(a)\equiv r_t\pmod{M_t}\}}\,\mathrm{d}a
 =\frac1{M_t}\int_I g(a)\,\mathrm{d}a
   +O_{I,t}(K A^{-3d/4}).
\]
Apply this identity successively to the largest remaining exponent. For the final exponent $m_1\ge d$, use \eqref{eq-prime-periodic-integral} directly. Dividing by $B-A$ gives
\begin{equation}\label{eq-prime-one-residue-choice}
 \mathbb P_I\bigl(u_{m_i}(a)\equiv r_i\pmod{M_i}
                  \text{ for every }i\bigr)
 =\frac1{M_1\cdots M_t}+O_{I,t}(K A^{-3d/4}).
\end{equation}
The same calculation directly proves \eqref{eq-prime-one-residue-choice} when $t=1$.

There are $M_1\cdots M_t\le\exp(Ct(\log K)^2)$ possible residue choices. Sum \eqref{eq-prime-one-residue-choice} over the choices in $F$. The resulting error is at most
\[
 O_{I,t}\bigl(K\exp(Ct(\log K)^2)A^{-3d/4}\bigr).
\]
Since $d\sim(\log K)^4$, this is at most $\exp(-c_I d)$ for all sufficiently large $K$, for example with $c_I=(\log A)/4$. This proves \eqref{eq-prime-joint-residue-error}.
\end{proof}

We next remove indices for which $u_m(a)$ has a small prime divisor. Put
\[
 z=\sqrt K,\qquad P=\prod_{p\le z}p,
 \qquad \rho=\prod_{p\le z}\left(1-\frac1p\right).
\]
These quantities depend on $K$. Mertens \cite{Mertens} proved the product estimate that gives
\begin{equation}\label{eq-prime-sieve-density}
 \rho=\frac{2e^{-\gamma}+o(1)}{\log K},
\end{equation}
where $\gamma$ is Euler's constant. For each exponent $m$, let
\[
 Y_m(a)=\mathbf1_{\{\gcd(u_m(a),P)=1\}}.
\]
The indices with $Y_m(a)=0$ will all be covered by the congruence $b\equiv0\pmod P$. The next lemma describes the joint distribution of the remaining indices and their residues modulo larger primes.

\begin{lemma}\label{lem-prime-residue-comparison}
Fix integers $t\ge1$ and $s\ge0$, and a real number $H>0$. Let $m_1,\ldots,m_t$ satisfy \eqref{eq-prime-separated-indices}, and let $q_1,\ldots,q_s$ be distinct primes in $[K,4K]$. For each pair $(i,j)$, let $X_{i,j}$ be a uniform random residue modulo $q_j$, with all these variables independent. Write $\mathbb P$ for probability with respect to these auxiliary variables. Let $F$ be any set of such residue arrays, and let $T\subseteq\{1,\ldots,t\}$. Then
\begin{equation}\label{eq-prime-residue-comparison}
 \begin{aligned}
 \Bigl|\mathbb P_I\bigl((u_{m_i}(a)\bmod q_j)_{i,j}\in F,\
                  Y_{m_i}(a)&=1\text{ for all }i\in T\bigr)\\
 &\qquad{}-\rho^{|T|}\mathbb P\bigl((X_{i,j})_{i,j}\in F\bigr)\Bigr|
 \le K^{-H}.
 \end{aligned}
\end{equation}
for all sufficiently large $K$, uniformly in all the choices above. When $s=0$, the event $F$ is either certain or impossible.
\end{lemma}

\begin{proof}
In this proof, we also write $F$ for the event that the residue array belongs to $F$. If $T$ is empty, the result follows from Lemma~\ref{lem-prime-joint-residues}, with the modulus $q_1\cdots q_s$ at every coordinate. Suppose now that $T$ is nonempty. The unwanted events are
\[
 p\mid u_{m_i}(a),\qquad i\in T,\quad p\le z.
\]
Let $R=2\lfloor(\log K)/2\rfloor$, so that $R$ is even and $R\sim\log K$. For $j\ge0$, let $S_j$ be the sum of the probabilities of the intersections of $F$ with $j$ distinct unwanted events. In particular, $S_0=\mathbb P_I(F)$. The finite inclusion and exclusion inequalities give
\begin{equation}\label{eq-prime-inclusion-exclusion}
 \sum_{j=0}^{R+1}(-1)^jS_j
 \le \mathbb P_I\bigl(F,\ Y_{m_i}=1\text{ for all }i\in T\bigr)
 \le \sum_{j=0}^{R}(-1)^jS_j.
\end{equation}
For completeness, at a point belonging to exactly $v\ge1$ unwanted events, the partial sum through order $h$ is
\[
 \sum_{j=0}^h(-1)^j\binom vj=(-1)^h\binom{v-1}{h}.
\]
This identity follows from Pascal's identity, with binomial coefficients beyond their range interpreted as zero. It is nonnegative when $h$ is even and nonpositive when $h$ is odd. At a point belonging to no unwanted event, both partial sums equal $1$. Multiplication by the indicator of $F$ and integration prove \eqref{eq-prime-inclusion-exclusion}.

Consider any term in these partial sums. At a given coordinate $m_i$, this term is determined by residues modulo an integer not exceeding
\[
 z^{R+1}(4K)^s=\exp\bigl(O_s((\log K)^2)\bigr).
\]
Indeed, at most $R+1$ small prime conditions occur in the term, and at most $s$ large primes are involved in $F$. Lemma~\ref{lem-prime-joint-residues} therefore compares this term, with error at most $\exp(-c_I d)$, to its probability when all coordinates are uniform modulo the indicated integers.

In this latter calculation, the residues at different coordinates are independent. At each coordinate, the residues modulo distinct primes are also independent by the Chinese remainder theorem. See Apostol \cite[Ch.~6]{Apostol}. We can therefore perform the calculation using independent variables that indicate divisibility by each small prime, each with probability $1/p$, together with the independent large prime residues $X_{i,j}$. In particular, the small prime variables are independent of $F$.

The total number of terms in the two partial sums is at most
\[
 (t\pi(z)+1)^{R+2}=\exp\bigl(O_t((\log K)^2)\bigr).
\]
Thus the total error from applying Lemma~\ref{lem-prime-joint-residues} to these terms is at most
\[
 \exp\bigl(-c_I d+O_t((\log K)^2)\bigr),
\]
which is smaller than any fixed inverse power of $K$.

It remains to bound the difference between the two partial sums in the independent calculation. Their difference is the sum of the probabilities of all intersections of $R+1$ distinct unwanted events, with the additional event $F$. On omitting $F$, this is at most
\[
 \frac1{(R+1)!}
 \left(t\sum_{p\le z}\frac1p\right)^{R+1}.
\]
The factor $(R+1)!$ is justified by expanding the power on the right. Every product over $R+1$ distinct unwanted events occurs $(R+1)!$ times, and the other products are nonnegative. Mertens \cite{Mertens} also proved
\[
 \sum_{p\le z}\frac1p=\log\log z+O(1).
\]
Using $(R+1)!\ge((R+1)/e)^{R+1}$, the difference is consequently at most
\[
 \left(\frac{O_t(\log\log K)}{R+1}\right)^{R+1}.
\]
The logarithm of this expression, divided by $\log K$, tends to $-\infty$. Hence this error too is smaller than any fixed inverse power of $K$.

In the independent calculation, the probability that none of the unwanted events occurs is exactly $\rho^{|T|}$, independently of $F$. Apply \eqref{eq-prime-inclusion-exclusion} to both calculations and use the two error estimates just obtained. Their sum is at most $K^{-H}$ for sufficiently large $K$. This proves \eqref{eq-prime-residue-comparison}.
\end{proof}

\subsection{A finite probability calculation}

A prime can cover several indices when their integer parts have the same residue modulo that prime. We shall choose many such groups so that no index and no prime is used twice. The finite calculation below gives the estimates needed for this choice.

Fix an integer $L\ge2$. Let $V$ be a set of $N$ indices, where $1\le N\le K$, and let $\mathcal R$ be a set of $r\le K$ primes in $[K,4K]$. For each $v\in V$, choose a random variable $Y_v$ that equals $1$ with probability $\rho$ and $0$ otherwise. Independently, for each $v\in V$ and $p\in\mathcal R$, choose a uniform residue $X_{v,p}$ modulo $p$. All these choices are independent. We write $\mathbb P$ and $\mathbb E$ for probability and expectation with respect to these finite independent choices and not to the choice of the real base $a$.

Put
\[
 T=\{v\in V\mid Y_v=1\}.
\]
We regard the elements of $T$ and the primes in $\mathcal R$ as belonging to two disjoint sets. For $p\in\mathcal R$ and $S\subseteq T$ with $|S|=L$, call $(p,S)$ admissible when
\[
 X_{v,p}=X_{w,p}\qquad(v,w\in S).
\]
We identify this choice with the $L+1$ element set consisting of $p$ and the elements of $S$. For a prime $p$, let $Z_p$ be the number of admissible choices that use $p$. Its expectation is
\begin{equation}\label{eq-prime-fixed-prime-mean}
 \nu_p=\binom NL\rho^L p^{-(L-1)}.
\end{equation}
Indeed, there are $\binom NL$ possible sets $S$. All their indices belong to $T$ with probability $\rho^L$, and their residues modulo $p$ are equal with probability $p^{-(L-1)}$.

For an index $v\in V$, first omit the requirement $Y_v=1$ and put
\[
 Z_v=\sum_{p\in\mathcal R}
      \sum_{\substack{S\subseteq V\setminus\{v\}\\|S|=L-1}}
      \prod_{i\in S}Y_i\mathbf1_{\{X_{i,p}=X_{v,p}\}}.
\]
When $Y_v=1$, this is the number of admissible choices that use $v$. Conditional on $Y_v=1$ and on all the residues $X_{v,p}$, its expectation is
\begin{equation}\label{eq-prime-index-mean}
 \mu=\binom{N-1}{L-1}\rho^{L-1}
       \sum_{p\in\mathcal R}p^{-(L-1)}.
\end{equation}
The right hand side does not depend on the residues on which we conditioned.

For two distinct objects $x,y$ belonging to the disjoint union of $V$ and $\mathcal R$, let $Z_{x,y}$ be the number of admissible choices containing both. We put $Z_{x,y}=0$ when one of the index objects involved does not belong to $T$.

\begin{lemma}\label{lem-prime-moments}
With the notation above, we have
\begin{align}\label{eq-prime-count-moments}
 \mathbb E(Z_p-\nu_p)^4&\ll_L K^2
          \qquad(p\in\mathcal R),\notag\\
 \mathbb E\bigl[Y_v(Z_v-\mu)^4\bigr]&\ll_L K^2
          \qquad(v\in V),
\end{align}
and
\begin{equation}\label{eq-prime-size-moment}
 \mathbb E(|T|-\rho N)^4\ll K^2.
\end{equation}
For any two distinct possible objects $x,y$, we also have
\begin{equation}\label{eq-prime-overlap-moment}
 \mathbb E Z_{x,y}^4\ll_L1.
\end{equation}
The constants are independent of $N,r,K,\rho$ and of the primes in $\mathcal R$.
\end{lemma}

\begin{proof}
If $\rho=0$, the assertions are immediate. We may therefore suppose that $\rho>0$. A fourth power of a sum expands into an ordered sum over four choices of summands. For fixed $L$, there are only finitely many possible patterns of coincidences among the index and prime symbols occurring in four summands. We fix one such pattern and count the possible assignments.

First fix $p\in\mathcal R$. For a set $S\subseteq V$ of size $L$, let $I_S$ indicate that every index of $S$ belongs to $T$ and that their residues modulo $p$ are equal. By \eqref{eq-prime-fixed-prime-mean},
\[
 Z_p-\nu_p
 =\sum_{\substack{S\subseteq V\\|S|=L}}
       (I_S-\mathbb E I_S).
\]
Consider a term corresponding to $S_1,S_2,S_3,S_4$. If one of these sets is disjoint from the other three, its centered factor is independent of the remaining factors and has mean zero. Such a term vanishes.

For a term that does not vanish, group the four sets by joining two groups whenever the corresponding sets intersect and repeating this operation until no further joining is possible. Every final group contains at least two of the four sets. Hence there are at most two groups. Let $c\le2$ be their number and let $m$ be the number of distinct indices in $S_1\cup\cdots\cup S_4$. If all four indicators equal $1$, then all these $m$ indices belong to $T$, and within each of the $c$ groups all residues modulo $p$ are equal. Therefore
\[
 \mathbb E\prod_{i=1}^4 I_{S_i}=\rho^m p^{-(m-c)}.
\]

We compare the centered product with this probability. Suppose some of the conditions $I_S=1$ have already been imposed. If a further set $S$ contains $b$ new indices and meets $h\ge1$ of the residue groups already present, its conditional probability is
\[
 \rho^b p^{-(b+h-1)}.
\]
Since $b\le L$ and $b+h\le L$, this is at least $\rho^L p^{-(L-1)}=\mathbb E I_S$. If $S$ is disjoint from the indices already present, its conditional probability is unchanged. Expanding the four centered factors therefore gives
\[
 \left|\mathbb E\prod_{i=1}^4(I_{S_i}-\mathbb E I_{S_i})\right|
 \le 2^4\mathbb E\prod_{i=1}^4I_{S_i}.
\]
For the fixed coincidence pattern there are at most $N^m$ assignments. Since $N\le K$, $p\ge K$ and $0\le\rho\le1$, its total contribution is at most
\[
 2^4N^m\rho^m p^{-(m-c)}\le2^4K^c\le2^4K^2.
\]
Summing over the finitely many patterns proves the first estimate in \eqref{eq-prime-count-moments}.

Now fix $v\in V$ and condition on $Y_v=1$ and all residues $X_{v,p}$. For $i\ne v$, the variables
\[
 Y_i,
 \qquad
 B_{i,p}=\mathbf1_{\{X_{i,p}=X_{v,p}\}}
 \qquad(p\in\mathcal R)
\]
are independent, with means $\rho$ and $1/p$. A summand of $Z_v$ has the form
\[
 I_{p,S}=\prod_{i\in S}Y_iB_{i,p},
 \qquad
 p\in\mathcal R,
 \qquad
 S\subseteq V\setminus\{v\},
 \qquad
 |S|=L-1.
\]
If the index set of one of four summands is disjoint from the index sets of the other three, its centered factor is independent of the remaining factors and the corresponding term vanishes.

Consider a term that need not vanish. Let $m$ be the number of distinct indices other than $v$ that occur, let $t$ be the number of distinct primes that occur, and let $q$ be the number of distinct pairs $(i,p)$ for which the condition $B_{i,p}=1$ occurs. Divide the four summands into classes by repeatedly joining two classes when the corresponding summands share an index or a prime. Since every summand shares an index with another summand, each final class contains at least two summands. There are therefore at most two classes.

Within one class, order the summands so that each one after the first shares an index or a prime with the preceding union. Such an ordering exists by the definition of the class. The first summand contains $L$ symbols and $L-1$ distinct pairs $(i,p)$. If a later summand introduces $u$ new symbols, it creates at least $u$ new pairs. Each new index is paired with its prime, and a new prime is paired with an old shared index. Thus a class containing $m_0$ distinct indices and $t_0$ distinct primes contains at least $m_0+t_0-1$ distinct pairs. Summing over the at most two classes gives
\[
 q\ge m+t-2.
\]
The probability that all four uncentered summands equal $1$ is at most $K^{-q}$. As above, conditioning some of the constituent independent variables to equal $1$ cannot decrease the probability of a further summand. The absolute value of the centered expectation is therefore at most $2^4K^{-q}$.

For the fixed coincidence pattern there are at most $N^m r^t\le K^{m+t}$ assignments. Its total contribution is consequently at most
\[
 2^4K^{m+t-q}\le2^4K^2.
\]
Summing over the finitely many patterns proves the second estimate in \eqref{eq-prime-count-moments} after removing the conditioning and multiplying by $\mathbb P(Y_v=1)\le1$.

For the size of $T$, write
\[
 |T|-\rho N=\sum_{v\in V}(Y_v-\rho).
\]
A term in the fourth moment vanishes if some index occurs exactly once. Every term that remains therefore contains at most two distinct indices. There are $O(N^2)$ such ordered terms, and each expectation has absolute value at most $1$. This proves \eqref{eq-prime-size-moment}.

It remains to prove \eqref{eq-prime-overlap-moment}. First suppose that $x$ and $y$ are two distinct indices. The count is zero unless both indices belong to $T$. Condition on this event and on all residues of the first index. A summand chooses a prime $p$, requires the second index to have the same residue modulo $p$, and chooses $L-2$ further retained indices with that residue. In a term of the fourth moment, let $t$ be the number of distinct chosen primes, let $m$ be the number of distinct additional indices, and let $q$ be the number of distinct pairs formed by one of these indices and one of the chosen primes at which a residue agreement is required. The probability is at most
\[
 K^{-t-q}.
\]
There are at most $N^m r^t\le K^{m+t}$ assignments for a fixed coincidence pattern. Since every additional index occurs in at least one of the $q$ pairs, we have $q\ge m$. Hence the contribution is at most $1$.

Next suppose that one object is an index $v$ and the other is a prime $p$. After conditioning on $Y_v=1$ and $X_{v,p}$, a term in the fourth moment involving $m$ distinct further indices has probability at most $p^{-m}\le K^{-m}$ and at most $N^m\le K^m$ assignments. Its contribution is again bounded. Two distinct primes never occur in the same admissible choice. Summing over the finitely many coincidence patterns proves \eqref{eq-prime-overlap-moment}.
\end{proof}

We shall use the following consequence of Pippenger and Spencer \cite{PippengerSpencer}. It is stated directly in the form needed below.

\begin{lemma}\label{lem-prime-disjoint-selection}
Fix an integer $k\ge2$ and a real number $0<\varepsilon<1$. There are $\delta>0$ and $D_0>0$ with the following property. Let $\Omega$ be a finite set and let $\mathcal F$ be a family of subsets of $\Omega$, each of cardinality $k$. Suppose that $D\ge D_0$, that every $x\in\Omega$ belongs to between $(1-\delta)D$ and $(1+\delta)D$ members of $\mathcal F$, and that every two distinct elements of $\Omega$ belong together to at most $\delta D$ members of $\mathcal F$. Then one can choose pairwise disjoint members of $\mathcal F$ whose union contains at least $(1-\varepsilon)|\Omega|$ elements.
\end{lemma}

\begin{proof}
Let $\theta>0$ be small enough that $(1+\theta)(1+\varepsilon/8)\le1+\varepsilon/2$. The theorem of Pippenger and Spencer \cite{PippengerSpencer} implies that, after decreasing $\delta$ and increasing $D_0$ if necessary, the family $\mathcal F$ can be partitioned into at most $(1+\theta)\Delta$ subfamilies, each consisting of pairwise disjoint sets, where $\Delta$ is the largest number of members of $\mathcal F$ containing one element of $\Omega$. We may also require $\delta\le\varepsilon/8$. Since $\Delta\le(1+\delta)D$, the number of these subfamilies is at most $(1+\varepsilon/2)D$.

Counting incidences between elements of $\Omega$ and members of $\mathcal F$ gives
\[
 k|\mathcal F|\ge(1-\delta)D|\Omega|.
\]
At least one of the subfamilies contains at least the average number of members. Since its members are pairwise disjoint, their union has at least
\[
 \frac{k|\mathcal F|}{(1+\varepsilon/2)D}
 \ge\frac{1-\delta}{1+\varepsilon/2}|\Omega|
 \ge(1-\varepsilon)|\Omega|
\]
elements. This proves the lemma.
\end{proof}

\subsection{Construction of the congruences}

We retain the interval $I=[A,B]$ and the quantities $d,z,P,\rho$ introduced above. We first keep the number $L$ of indices covered by one prime fixed.

\begin{proposition}\label{prop-prime-fixed-L}
Fix an integer $L\ge2$. For almost every $a\in I$ and every sufficiently large integer $K$ that is a power of two, there are a squarefree integer $Q\ge2$ and an integer $b$ such that all prime divisors of $Q$ are at most $8K$,
\[
 \gcd\bigl(b-u_m(a),Q\bigr)>1\qquad(1\le m\le K),
\]
and
\begin{equation}\label{eq-prime-fixed-L-size}
 \log Q\le\left(\frac{2e^{-\gamma}}L+o(1)\right)K.
\end{equation}
The term $o(1)$ may depend on $a,I,L$, and the limit is taken through powers of two.
\end{proposition}

\begin{proof}
First set aside the indices $1,\ldots,d-1$. Partition the remaining indices $\{d,\ldots,K\}$ into sets $V_1,\ldots,V_d$ according to their residues modulo $d$, and put $N_j=|V_j|$. Then
\[
 N_j=(1+o(1))\frac Kd
\]
uniformly for $1\le j\le d$. The least index in each set is at least $d$, and distinct indices in the same set differ by at least $d$. Thus Lemma~\ref{lem-prime-residue-comparison} applies to any fixed number of indices in one set $V_j$.

Allocate
\[
 r_j=\left\lfloor\frac{\rho N_j}{L}\right\rfloor
\]
primes to $V_j$. For all sufficiently large $K$, every $r_j$ is positive because $\rho N_j\asymp K/(\log K)^5\to\infty$. By \eqref{eq-prime-sieve-density}, the total number of primes needed is
\begin{equation}\label{eq-prime-number-used}
 R=\sum_{j=1}^d r_j
  =\frac\rho L(K-d+1)+O(d)
  =\left(\frac{2e^{-\gamma}}L+o(1)\right)\frac K{\log K}.
\end{equation}
Choose the first $R$ primes not smaller than $K$, and divide this consecutive list into disjoint sets $\mathcal R_j$ of sizes $r_j$, in that order. The prime number theorem \cite[Ch.~14]{Apostol} gives
\[
 \pi(4K)-\pi(K)=(3+o(1))\frac K{\log K}.
\]
Since $2e^{-\gamma}/L<1$ for $L\ge2$, \eqref{eq-prime-number-used} shows that all the chosen primes lie in $[K,4K]$ for large $K$.

Let $q_j$ be the smallest prime in $\mathcal R_j$. We claim that
\begin{equation}\label{eq-prime-ranges}
 \max_{p\in\mathcal R_j}\left|\frac p{q_j}-1\right|=o(1)
 \quad\text{uniformly in }j.
\end{equation}
If $p_n$ denotes the $n$th prime, the prime number theorem gives $p_n\sim n\log n$. All prime indices in the chosen list are between positive constant multiples of $K/\log K$, whereas
\[
 \max_j r_j=O_L\left(\frac K{d\log K}\right)
            =o(K/\log K).
\]
Hence the first and last prime indices in every $\mathcal R_j$ have ratio tending to $1$, uniformly in $j$. The asymptotic formula for $p_n$ now proves \eqref{eq-prime-ranges}.

For $a\in I$, put
\[
 T_j(a)=\{m\in V_j\mid Y_m(a)=1\}.
\]
For $p\in\mathcal R_j$ and $S\subseteq T_j(a)$ with $|S|=L$, call $(p,S)$ admissible if the integers $u_m(a)$, $m\in S$, have the same residue modulo $p$. We regard the indices in $T_j(a)$ and the primes in $\mathcal R_j$ as elements of two disjoint sets, and identify $(p,S)$ with the $L+1$ element set consisting of these objects.

Use the deterministic quantities in \eqref{eq-prime-fixed-prime-mean} and \eqref{eq-prime-index-mean} with $N=N_j$ and $\mathcal R=\mathcal R_j$, and denote them by $\nu_p$ and $\mu_j$. Put
\[
 D_j=\binom{N_j}{L}\rho^L q_j^{-(L-1)}.
\]
Equations \eqref{eq-prime-ranges} and \eqref{eq-prime-fixed-prime-mean} give $\nu_p=(1+o(1))D_j$ uniformly for $p\in\mathcal R_j$. By \eqref{eq-prime-index-mean},
\[
 \frac{\mu_j}{D_j}
 =\frac{L}{\rho N_j}
   \sum_{p\in\mathcal R_j}\left(\frac{q_j}{p}\right)^{L-1}
 =(1+o(1))\frac{Lr_j}{\rho N_j}
 =1+o(1).
\]
All the errors are uniform in $j$. Since $N_j\asymp K/d$, $q_j\in[K,4K]$ and \eqref{eq-prime-sieve-density} holds, we have
\begin{equation}\label{eq-prime-common-counts}
 \nu_p=(1+o(1))D_j,
 \qquad
 \mu_j=(1+o(1))D_j,
 \qquad
 D_j\asymp_L\frac K{(\log K)^{5L}}.
\end{equation}
In particular, the smallest $D_j$ tends to infinity.

We next transfer the fourth moment estimates of Lemma~\ref{lem-prime-moments} to the actual powers. For an object $x\in T_j(a)\sqcup\mathcal R_j$, let $Z_j(a,x)$ be the number of admissible choices containing $x$. For two distinct possible objects $x,y\in V_j\sqcup\mathcal R_j$, let $Z_j(a,x,y)$ be the number containing both, and put this number equal to zero whenever an index among $x,y$ does not belong to $T_j(a)$.

A count attached to one object is a sum of at most $K^L$ indicators, and its mean in the independent calculation is at most $K^L$. If $Z$ is such a count with independent mean $\lambda$, then
\[
 (Z-\lambda)^4=\sum_{h=0}^4\binom4h(-\lambda)^{4-h}Z^h.
\]
After expanding the powers $Z^h$, the sum of the absolute values of all coefficients is at most
\[
 \sum_{h=0}^4\binom4hK^{L(4-h)}K^{Lh}
 =2^4K^{4L}.
\]
The same bound is sufficient for the fourth moment of $|T_j(a)|-\rho N_j$ and for the ordinary fourth moments of the pair counts. Every resulting term involves at most $4L+2$ indices from $V_j$ and at most four primes from $\mathcal R_j$. Repeated indicators may be replaced by one indicator.

Apply \eqref{eq-prime-residue-comparison} to every resulting term with $H=4L+10$. The difference between its expectation over $a\in I$ and its expectation in the independent calculation is at most $K^{-4L-10}$. Hence the total difference for any expanded fourth moment is at most $2^4K^{-10}$. Lemma~\ref{lem-prime-moments} therefore gives the same fourth moment bounds for the actual powers, after changing the implied constants. For an index $v\in V_j$, put
\[
 Z_v(a)=\sum_{p\in\mathcal R_j}
      \sum_{\substack{S\subseteq V_j\setminus\{v\}\\|S|=L-1}}
      \prod_{i\in S}Y_i(a)
      \mathbf1_{\{u_i(a)\equiv u_v(a)\pmod p\}}.
\]
The relevant estimate is applied to
\[
 \mathbb E_I\bigl[Y_v(a)(Z_v(a)-\mu_j)^4\bigr].
\]
On the event $Y_v(a)=1$, the quantity $Z_v(a)$ is exactly $Z_j(a,v)$.

For a random variable $X$ and $t>0$,
\[
 \mathbb P(|X|>t)\le t^{-4}\mathbb E|X|^4.
\]
By \eqref{eq-prime-common-counts}, the threshold $D_j/\log K$ is at least a positive constant, depending only on $L$, times $K/(\log K)^{5L+1}$. Thus \eqref{eq-prime-count-moments} shows that the probability of failure for any one prime or index is
\[
 O_{I,L}\left(K^{-2}(\log K)^{20L+4}\right).
\]
There are at most $K+R\le2K$ such objects over all $j$ for large $K$.

Also, $\rho N_j/\log K\gg_I K/(\log K)^6$. By \eqref{eq-prime-size-moment}, the probability that the required estimate for $|T_j(a)|$ fails is $O_I(K^{-2}(\log K)^{24})$. There are $d=O((\log K)^4)$ sets $V_j$. Finally, \eqref{eq-prime-overlap-moment} shows that the probability of
\[
 Z_j(a,x,y)>K^{3/4}
\]
for one fixed pair is $O_{I,L}(K^{-3})$. There are at most $(K+R)^2\le4K^2$ pairs over all $j$.

It follows that, except on a set of bases of probability
\begin{equation}\label{eq-prime-bad-probability}
 O_{I,L}\left(K^{-1}(\log K)^{20L+4}\right),
\end{equation}
the following estimates hold simultaneously for every $j$,
\begin{align}\label{eq-prime-simultaneous-estimates}
 \bigl||T_j(a)|-\rho N_j\bigr|&\le\frac{\rho N_j}{\log K},\notag\\
 |Z_j(a,p)-\nu_p|&\le\frac{D_j}{\log K}
                 \qquad(p\in\mathcal R_j),\notag\\
 |Z_j(a,v)-\mu_j|&\le\frac{D_j}{\log K}
                 \qquad(v\in T_j(a)),\notag\\
 Z_j(a,x,y)&\le K^{3/4}
                 \qquad(x\ne y).
\end{align}
We now restrict $K$ to powers of two. The series obtained by summing \eqref{eq-prime-bad-probability} over dyadic $K$ converges. Since $\mathbb P_I$ is normalized Lebesgue measure on $I$, Lemma~\ref{lem-borel-cantelli} shows that, for almost every fixed $a\in I$, all the estimates in \eqref{eq-prime-simultaneous-estimates} hold for every sufficiently large dyadic $K$.

Fix one such base $a$. By \eqref{eq-prime-common-counts} and \eqref{eq-prime-simultaneous-estimates}, every object in $T_j(a)\sqcup\mathcal R_j$ belongs to $(1+o(1))D_j$ admissible choices, uniformly in $j$. Any two distinct objects belong together to at most $K^{3/4}$ admissible choices, and
\[
 \frac{K^{3/4}}{D_j}
 \ll_L\frac{(\log K)^{5L}}{K^{1/4}}
 \longrightarrow0
\]
uniformly in $j$. Lemma~\ref{lem-prime-disjoint-selection}, with $k=L+1$, now gives a collection of pairwise disjoint admissible choices for each $j$ which leaves only
\[
 o(|T_j(a)|+r_j)
\]
objects unused. To see the uniformity, fix $\varepsilon>0$ in that lemma. For all sufficiently large dyadic $K$, its hypotheses hold for every $j$ with the same $\delta$ and $D_0$. The number of unused objects is then at most $\varepsilon(|T_j(a)|+r_j)$ for every $j$. Since $\varepsilon$ is arbitrary, the stated estimate follows uniformly.

Consequently, the total number of retained indices which are not contained in one of the selected admissible choices is
\begin{equation}\label{eq-prime-uncovered-indices}
 o\left(\sum_{j=1}^d(|T_j(a)|+r_j)\right)
 =o(\rho K)=o(K/\log K).
\end{equation}
Here we used \eqref{eq-prime-simultaneous-estimates} and \eqref{eq-prime-number-used}. Since the selected admissible choices are pairwise disjoint and the sets $\mathcal R_j$ are disjoint, no prime is assigned to two different groups of indices.

We can now prescribe the congruences for $b$. For every prime $p\le z$, require $b\equiv0\pmod p$. This covers every index $m$ with $Y_m(a)=0$. For every selected admissible choice $(p,S)$, require $b$ to equal the common residue of the integers $u_m(a)$, $m\in S$, modulo $p$. The remaining indices are the uncovered retained indices and the initially set aside indices $1,\ldots,d-1$. By \eqref{eq-prime-uncovered-indices}, their number is
\[
 o(K/\log K)+O((\log K)^4)=o(K/\log K).
\]
Assign a distinct prime in $(4K,8K]$ to each remaining index $m$, and require $b\equiv u_m(a)$ modulo its assigned prime. There are enough such primes because the prime number theorem gives
\[
 \pi(8K)-\pi(4K)=(4+o(1))\frac K{\log K}.
\]
All primes used in this step are different from the small primes and from the primes in the sets $\mathcal R_j$.

Let $Q$ be the product of all primes to which a residue has been assigned. It is squarefree, and all its prime divisors are at most $8K$. By the Chinese remainder theorem \cite[Ch.~6]{Apostol}, there is an integer $b$ satisfying all the prescribed congruences.

Every index $m\le K$ is now covered by at least one prime dividing both $Q$ and $b-u_m(a)$. It remains to estimate $Q$. The small primes contribute at most $\sqrt K\log K=o(K)$ to $\log Q$. The primes used by the selected admissible choices contribute at most $R\log(4K)$. The last group consists of $o(K/\log K)$ primes, each at most $8K$, and contributes $o(K)$. By \eqref{eq-prime-number-used},
\[
 \log Q
 \le \sqrt K\log K+R\log(4K)
                +o(K/\log K)\log(8K)
 \le\left(\frac{2e^{-\gamma}}L+o(1)\right)K.
\]
This proves \eqref{eq-prime-fixed-L-size}.
\end{proof}

\begin{proof}[Proof of Proposition~\ref{prop-prime-cover}]
For each integer $L\ge2$, Proposition~\ref{prop-prime-fixed-L} gives a subset of $I$ whose complement has measure zero. Intersect these sets over all $L\ge2$. The intersection still has full measure.

Fix a base $a$ in this intersection. For every sufficiently large dyadic integer $K$, consider the squarefree integers whose prime divisors are at most $8K$ and for which some residue $b$ satisfies \eqref{eq-prime-cover}. There are only finitely many such squarefree integers. Proposition~\ref{prop-prime-fixed-L} with $L=2$ shows that the collection is nonempty. Let $Q_K$ be its smallest member, and choose one corresponding residue $b_K$.

For every fixed $L\ge2$, minimality and \eqref{eq-prime-fixed-L-size} give
\[
 \limsup_{\substack{K\to\infty\\K\text{ a power of }2}}
 \frac{\log Q_K}{K}
 \le\frac{2e^{-\gamma}}L.
\]
The left hand side is nonnegative. Letting $L$ tend to infinity proves \eqref{eq-prime-cover-size}. The order of the two limits is important. We first fix $L$ and let the dyadic integer $K$ tend to infinity. No estimate with a growing $L$ is used.

This proves the proposition for almost every base in the fixed interval $I$. Finally, the intervals $[1+1/j,j]$, $j\ge2$, cover $(1,\infty)$. Applying the result on each of these intervals and taking the union of their exceptional sets proves the proposition for almost every $a>1$.
\end{proof}

The modulus in Proposition~\ref{prop-prime-cover} need not remain bounded. Accordingly, \eqref{eq-prime-progression-count} gives the growth estimate in Theorem~\ref{thm-prime-exceptions}, but does not imply the positive lower density required in Conjecture~\ref{conj-near-one-obstruction}.

\section{A quantitative squarefree complement}\label{sec-squarefree-proof}

    We first record the estimate for growing square moduli that will be used in the proof of \eqref{eq-squarefree-quantitative}.

    \begin{lemma}\label{lem-growing-squares}
        Let $J=[A_1,A_2]$, where $1<A_1<A_2$, and let $\varepsilon>0$. Put
        $$
        K_j=2^j,
        \qquad
        z_j=\frac{\sqrt{K_j}}{(\log K_j)^{1+\varepsilon}}.
        $$
        For a prime $p\le z_j$ and $0\le r<p^2$, let
        $$
        N_{j,p,r}(a)
        =\#\{m\le K_j\mid \lfloor a^m\rfloor\equiv r\pmod{p^2}\}.
        $$
        Then, for almost every $a\in J$,
        \begin{equation}\label{eq-growing-square-discrepancy}
            \sum_{p\le z_j}
            \max_{0\le r<p^2}
            \left|N_{j,p,r}(a)-\frac{K_j}{p^2}\right|
            =o_{a,\varepsilon}(K_j).
        \end{equation}
    \end{lemma}

    \begin{proof}
        For $p\le z_j$, apply Lemma~\ref{lem-erdos-turan} with $H=K_j$ to the points $a^m/p^2$ and the intervals
        $$
        \left[\frac r{p^2},\frac{r+1}{p^2}\right).
        $$
        As in the proof of Lemma~\ref{lem-residue-ud}, membership in this interval is equivalent to $\lfloor a^m\rfloor\equiv r\pmod{p^2}$. Since $p^2\le z_j^2<K_j$ for all sufficiently large $j$, Lemma~\ref{lem-exp-sum} with $B=1$ gives, uniformly in $p\le z_j$ and $1\le h\le K_j$,
        $$
        \int_J
        \left|\sum_{m\le K_j}e\left(\frac{h a^m}{p^2}\right)\right|^2\,\mathrm{d}a
        \ll_J K_j.
        $$
        The right hand side of Lemma~\ref{lem-erdos-turan} is independent of $r$. Taking the maximum over $r$, integrating, and using Cauchy's inequality therefore gives
        $$
        \int_J
        \max_{0\le r<p^2}
        \left|N_{j,p,r}(a)-\frac{K_j}{p^2}\right|\,\mathrm{d}a
        \ll_J \sqrt{K_j}\log K_j.
        $$
        Since $\log z_j\asymp\log K_j$, summing over $p\le z_j$ and using $\pi(z)\ll z/\log z$ gives
        \begin{equation}\label{eq-growing-square-L1}
            \int_J
            \sum_{p\le z_j}
            \max_{0\le r<p^2}
            \left|N_{j,p,r}(a)-\frac{K_j}{p^2}\right|\,\mathrm{d}a
            \ll_{J,\varepsilon}
            \frac{K_j}{(\log K_j)^{1+\varepsilon}}.
        \end{equation}
        For each positive integer $s$, Markov's inequality and \eqref{eq-growing-square-L1} give
        $$
        \sum_{j=1}^{\infty}
        \operatorname{meas}\left\{a\in J\mid
        \sum_{p\le z_j}
        \max_{0\le r<p^2}
        \left|N_{j,p,r}(a)-\frac{K_j}{p^2}\right|>
        \frac{K_j}{s}\right\}<\infty,
        $$
        since $\log K_j\asymp j$. Thus the exceptional measures are summable. Lemma~\ref{lem-borel-cantelli} and a countable intersection over $s$ now give \eqref{eq-growing-square-discrepancy}.\qedhere
    \end{proof}

    \begin{proof}[Proof of \eqref{eq-squarefree-quantitative}]
        Fix $\varepsilon>0$ and a compact interval $J=[A_1,A_2]\subset(1,\infty)$. Choose $a\in J$ for which Lemma~\ref{lem-growing-squares} holds, and put
        $$
        t_m=\lfloor a^m\rfloor.
        $$
        We shall use the elementary estimate
        \begin{equation}\label{eq-prime-square-sum}
            \sum_{p\in\mathcal{P}}\frac1{p^2}
            \le
            \frac14+
            \sum_{\substack{n\ge3\\n\ \mathrm{odd}}}\frac1{n^2}
            =
            \frac{\pi^2}{8}-\frac34
            <\frac12.
        \end{equation}

        Let $K_j$ and $z_j$ be as in Lemma~\ref{lem-growing-squares}. For every integer $n$, the number of indices $m\le K_j$ for which
        $$
        p^2\mid n-t_m
        $$
        for at least one prime $p\le z_j$ is at most
        $$
        \sum_{p\le z_j}N_{j,p,n\bmod p^2}(a)
        \le
        K_j\sum_{p\le z_j}\frac1{p^2}
        +
        \sum_{p\le z_j}
        \max_{0\le r<p^2}
        \left|N_{j,p,r}(a)-\frac{K_j}{p^2}\right|.
        $$
        It follows from \eqref{eq-growing-square-discrepancy} and \eqref{eq-prime-square-sum} that, for all sufficiently large $j$, uniformly in $n$, at least $K_j/2$ indices $m\le K_j$ satisfy
        \begin{equation}\label{eq-no-small-square}
            p^2\nmid n-t_m
            \qquad
            \text{for every prime }p\le z_j.
        \end{equation}

        Let $F_a(x)$ denote the number of positive integers not exceeding $x$ which do not belong to
        $$
        \mathcal{Q}+\{t_m\mid m\ge1\}.
        $$
        Put
        $$
        U_j=\lfloor a^{K_j}\rfloor.
        $$
        Since $U_j\to\infty$ and the sequence $U_j$ is eventually strictly increasing, for sufficiently large $x$ we may choose $j$ so that
        $$
        U_j\le x<U_{j+1}.
        $$
        Consider an exceptional integer $n$ with
        $$
        U_{j-1}<n\le x.
        $$
        Apply \eqref{eq-no-small-square} with $K_{j-1}$. For at least $K_{j-1}/2$ indices $m\le K_{j-1}$, the positive integer $n-t_m$ is not divisible by $p^2$ for any prime $p\le z_{j-1}$. Since $n$ is exceptional, none of these differences is squarefree. Hence each of them is divisible by $p^2$ for some prime $p>z_{j-1}$.

        For a fixed $m\le K_{j-1}$, the number of integers $n\le x$ for which
        $$
        p^2\mid n-t_m
        $$
        for some prime $p>z_{j-1}$ is at most
        $$
        \sum_{z_{j-1}<p\le\sqrt{x}}
        \left(\frac{x}{p^2}+1\right)
        \le
        \frac{x}{z_{j-1}-1}+\sqrt{x}.
        $$
        There are $F_a(x)-F_a(U_{j-1})$ exceptional integers in $(U_{j-1},x]$, and $F_a(U_{j-1})\le U_{j-1}$. Thus, whenever $F_a(x)>U_{j-1}$, their number is at least $F_a(x)-U_{j-1}$. Counting pairs $(n,m)$ therefore gives
        $$
        \frac{K_{j-1}}2
        \bigl(F_a(x)-U_{j-1}\bigr)
        \le
        K_{j-1}
        \left(\frac{x}{z_{j-1}-1}+\sqrt{x}\right).
        $$
        Thus, in all cases,
        \begin{equation}\label{eq-squarefree-dyadic}
            F_a(x)
            \le
            U_{j-1}
            +\frac{2x}{z_{j-1}-1}
            +2\sqrt{x}.
        \end{equation}

        Since $K_j=2^j$ and $U_j=\lfloor a^{K_j}\rfloor$, the choice of $j$ gives
        $$
        K_{j-1}\asymp_a\log x.
        $$
        Consequently
        $$
        z_{j-1}
        \asymp_{a,\varepsilon}
        \frac{\sqrt{\log x}}{(\log\log x)^{1+\varepsilon}}.
        $$
        Moreover, for all sufficiently large $j$,
        $$
        U_j\ge \frac12a^{K_j},
        $$
        so
        $$
        \frac{U_{j-1}}x
        \le \frac{U_{j-1}}{U_j}
        \ll_a a^{-K_{j-1}}.
        $$
        The term $x^{-1/2}$ is also negligible compared with the reciprocal of $z_{j-1}$. Hence \eqref{eq-squarefree-dyadic} gives
        $$
        F_a(x)
        \ll_{a,\varepsilon}
        \frac{x(\log\log x)^{1+\varepsilon}}{\sqrt{\log x}}.
        $$

        For each fixed $\varepsilon>0$, the conclusion holds for almost all $a$ in every compact subinterval of $(1,\infty)$. Take the countable intersection corresponding to the intervals $[1+1/n,n]$ and to $\varepsilon=1/s$, where $n,s\ge2$. If $\varepsilon>0$ is arbitrary, choose $s$ with $1/s<\varepsilon$. The estimate with $1/s$ implies \eqref{eq-squarefree-quantitative}. This proves \eqref{eq-squarefree-quantitative}.
    \end{proof}

    \subsection{Two integer parts of powers}

    For the second estimate, we need an exponential sum bound with two parameters. We separate the two ranges of exponents so that the larger exponent controls every off diagonal term.

    \begin{lemma}\label{lem-two-block-exp-sum}
        Let $J=[A_1,A_2]$, where $1<A_1<A_2$. There exists an integer $M=M(J)\ge2$ such that, uniformly for $K\ge2$ and $1\le d,h\le K^2$,
        $$
        \int_J
        \left|
        \sum_{m\le K}
        \sum_{MK<\ell\le(M+1)K}
        e\left(\frac{h(a^m+a^\ell)}d\right)
        \right|^2\,\mathrm{d}a
        \ll_J K^2.
        $$
    \end{lemma}

    \begin{proof}
        Choose $M$ so large that $A_1^M>A_2^2$. We may assume that $K$ is sufficiently large, since the remaining values can be absorbed into the implied constant. If $m_1,m_2\le K$ and
        $$
        MK<\ell_1<\ell_2\le(M+1)K,
        $$
        put
        $$
        F(a)=a^{m_1}-a^{m_2}+a^{\ell_1}-a^{\ell_2}.
        $$
        For $a\in J$,
        $$
        \ell_2a^{\ell_2-1}-\ell_1a^{\ell_1-1}
        \ge (\ell_2-\ell_1)a^{\ell_2-1}
        \ge A_1^{\ell_2-1},
        $$
        while
        $$
        \left|m_1a^{m_1-1}-m_2a^{m_2-1}\right|
        \le 2KA_2^{K-1}.
        $$
        Similarly,
        $$
        \begin{aligned}
        \ell_2(\ell_2-1)a^{\ell_2-2}
        -\ell_1(\ell_1-1)a^{\ell_1-2}
        \ge
        \bigl(\ell_2(\ell_2-1)-\ell_1(\ell_1-1)\bigr)
        A_1^{\ell_2-2}
        \ge 2MKA_1^{\ell_2-2},
        \end{aligned}
        $$
        and the absolute value of the second derivative of $a^{m_1}-a^{m_2}$ is at most $2K^2A_2^{K-2}$. It follows from the choice of $M$ that, for all sufficiently large $K$,
        $$
        |F'(a)|\ge \frac12A_1^{\ell_2-1}
        $$
        and $F''(a)<0$ on $J$. Thus $F'$ is monotone. The first derivative estimate therefore gives
        \begin{equation}\label{eq-two-block-off-diagonal}
        \left|
        \int_J e\left(\frac{hF(a)}d\right)\,\mathrm{d}a
        \right|
        \ll_J \frac{d}{hA_1^{\ell_2-1}}.
        \end{equation}
        The same bound holds when $\ell_1>\ell_2$, with $\ell_2$ replaced by $\ell_1$.

        We now expand the square. The terms with $\ell_1=\ell_2$ contribute
        $$
        K\int_J
        \left|\sum_{m\le K}e\left(\frac{ha^m}d\right)\right|^2\,\mathrm{d}a
        \ll_J K^2.
        $$
        Indeed, the proof of Lemma~\ref{lem-exp-sum} gives an upper bound $O_J(K+1+\log^2(2d/h))$, which is $O_J(K)$ under $d,h\le K^2$. By \eqref{eq-two-block-off-diagonal}, the terms with $\ell_1\ne\ell_2$ contribute at most
        $$
        \frac dh K^2
        \sum_{\substack{MK<\ell_1,\ell_2\le(M+1)K\\
        \ell_1\ne\ell_2}}
        A_1^{-\max(\ell_1,\ell_2)+1}
        \ll_J K^5A_1^{-MK}
        \ll_J1.
        $$
        This proves the lemma.
    \end{proof}

    \begin{lemma}\label{lem-growing-square-pairs}
        Let $J=[A_1,A_2]$, where $1<A_1<A_2$, let $M$ be as in Lemma~\ref{lem-two-block-exp-sum}, and let $\varepsilon>0$. Put
        $$
        K_j=2^j,
        \qquad
        z_j=\frac{K_j}{(\log K_j)^{1+\varepsilon}}.
        $$
        For a prime $p\le z_j$ and $0\le r<p^2$, let
        $$
        R_{j,p,r}(a)
        =
        \#\left\{(m,\ell)\mid
        \begin{array}{c}
        1\le m\le K_j,\quad MK_j<\ell\le(M+1)K_j,\\
        \lfloor a^m+a^\ell\rfloor\equiv r\pmod{p^2}
        \end{array}
        \right\}.
        $$
        Then, for almost every $a\in J$,
        \begin{equation}\label{eq-growing-square-pair-discrepancy}
        \sum_{p\le z_j}
        \max_{0\le r<p^2}
        \left|R_{j,p,r}(a)-\frac{K_j^2}{p^2}\right|
        =o_{a,\varepsilon}(K_j^2).
        \end{equation}
    \end{lemma}

    \begin{proof}
        For $p\le z_j$, apply Lemma~\ref{lem-erdos-turan} with $H=K_j^2$ to the $K_j^2$ points
        $$
        \frac{a^m+a^\ell}{p^2},
        \qquad
        1\le m\le K_j,
        \quad
        MK_j<\ell\le(M+1)K_j.
        $$
        Since $p^2<K_j^2$ for all sufficiently large $j$, Lemma~\ref{lem-two-block-exp-sum} and Cauchy's inequality give
        $$
        \int_J
        \max_{0\le r<p^2}
        \left|R_{j,p,r}(a)-\frac{K_j^2}{p^2}\right|\,\mathrm{d}a
        \ll_J K_j\log K_j.
        $$
        Consequently,
        \begin{align}\label{eq-growing-square-pair-L1}
        &\int_J
        \sum_{p\le z_j}
        \max_{0\le r<p^2}
        \left|R_{j,p,r}(a)-\frac{K_j^2}{p^2}\right|\,\mathrm{d}a
        \ll_{J,\varepsilon}
        K_j\log K_j\,\pi(z_j)
        \ll_{J,\varepsilon}
        \frac{K_j^2}{(\log K_j)^{1+\varepsilon}}.
        \end{align}
        For each positive integer $s$, Markov's inequality and \eqref{eq-growing-square-pair-L1} show that the measures of the sets on which the sum in \eqref{eq-growing-square-pair-discrepancy} exceeds $K_j^2/s$ are summable over $j$. Lemma~\ref{lem-borel-cantelli} and a countable intersection over $s$ prove the assertion.
    \end{proof}

    \begin{proof}[Proof of \eqref{eq-squarefree-two-powers}]
        Fix $\varepsilon>0$ and a compact interval $J=[A_1,A_2]\subset(1,\infty)$. Let $M=M(J)$ be as in Lemma~\ref{lem-two-block-exp-sum}, and choose $a\in J$ for which Lemma~\ref{lem-growing-square-pairs} holds. Put
        $$
        t_m=\lfloor a^m\rfloor,
        \qquad
        s_{m,\ell}=\lfloor a^m+a^\ell\rfloor.
        $$
        Then
        \begin{equation}\label{eq-two-floor-alternatives}
            t_m+t_\ell\in\{s_{m,\ell}-1,s_{m,\ell}\}.
        \end{equation}
        By \eqref{eq-prime-square-sum},
        $$
        \eta=1-2\sum_{p\in\mathcal{P}}\frac1{p^2}>0.
        $$
        Let $K_j$ and $z_j$ be as in Lemma~\ref{lem-growing-square-pairs}. For every integer $n$, the number of pairs $(m,\ell)$ in the ranges defining $R_{j,p,r}(a)$ for which
        $$
        p^2\mid n-s_{m,\ell}
        \quad\text{or}\quad
        p^2\mid n-s_{m,\ell}+1
        $$
        for at least one prime $p\le z_j$ is at most
        $$
        2K_j^2\sum_{p\le z_j}\frac1{p^2}
        +2\sum_{p\le z_j}
        \max_{0\le r<p^2}
        \left|R_{j,p,r}(a)-\frac{K_j^2}{p^2}\right|.
        $$
        The consecutive integers $n-s_{m,\ell}$ and $n-s_{m,\ell}+1$ are relatively prime. A prime square divides their product only if it divides one of the two factors. It follows from \eqref{eq-growing-square-pair-discrepancy} that, for all sufficiently large $j$, uniformly in $n$, at least $\eta K_j^2/2$ pairs satisfy
        \begin{equation}\label{eq-two-no-small-squares}
            p^2\nmid
        \bigl(n-s_{m,\ell}\bigr)
        \bigl(n-s_{m,\ell}+1\bigr)
        \qquad
        \text{for every prime }p\le z_j.
        \end{equation}

        Let $F_{a,2}(x)$ denote the number of positive integers not exceeding $x$ which do not belong to
        $$
        \mathcal{Q}+\{t_m+t_\ell\mid 1\le m<\ell\}.
        $$
        Put
        $$
        V_j=\left\lceil2a^{(M+1)K_j}\right\rceil.
        $$
        Since $V_j\to\infty$ and the sequence $V_j$ is eventually strictly increasing, for sufficiently large $x$ choose $j$ so that
        $$
        V_j\le x<V_{j+1}.
        $$
        Consider an exceptional integer $n$ with
        $$
        V_{j-1}<n\le x.
        $$
        Apply \eqref{eq-two-no-small-squares} with $K_{j-1}$. For every pair counted there, we have $s_{m,\ell}\le V_{j-1}<n$. Put
        $$
        q=n-s_{m,\ell}.
        $$
        If both $q$ and $q+1$ were squarefree, then \eqref{eq-two-floor-alternatives} would give a representation of $n$ as a squarefree integer plus $t_m+t_\ell$. Hence, for every one of the at least $\eta K_{j-1}^2/2$ pairs in \eqref{eq-two-no-small-squares}, there is a prime $p>z_{j-1}$ such that
        $$
        p^2\mid q
        \quad\text{or}\quad
        p^2\mid q+1.
        $$
        For a fixed pair $(m,\ell)$, the number of integers $n\le x$ for which this happens is at most
        $$
        2\sum_{z_{j-1}<p\le\sqrt{x+1}}
        \left(\frac{x}{p^2}+1\right)
        \ll
        \frac{x}{z_{j-1}-1}+\sqrt{x}.
        $$
        There are $F_{a,2}(x)-F_{a,2}(V_{j-1})$ exceptional integers in $(V_{j-1},x]$, and $F_{a,2}(V_{j-1})\le V_{j-1}$. Counting pairs $(n,(m,\ell))$ therefore gives
        $$
        F_{a,2}(x)
        \le
        V_{j-1}
        +O_{a,J}\left(\frac{x}{z_{j-1}}+\sqrt{x}\right).
        $$

        The choice of $j$ gives
        $$
        K_{j-1}\asymp_{a,J}\log x
        $$
        and hence
        $$
        z_{j-1}
        \asymp_{a,J,\varepsilon}
        \frac{\log x}{(\log\log x)^{1+\varepsilon}}.
        $$
        Moreover,
        $$
        \frac{V_{j-1}}x
        \le\frac{V_{j-1}}{V_j}
        \ll_{a,J}a^{-(M+1)K_{j-1}},
        $$
        and the term $x^{-1/2}$ is negligible compared with $1/z_{j-1}$. We conclude that
        $$
        F_{a,2}(x)
        \ll_{a,\varepsilon}
        \frac{x(\log\log x)^{1+\varepsilon}}{\log x}.
        $$

        For each fixed $\varepsilon>0$, this holds for almost all $a$ in every compact subinterval of $(1,\infty)$. Taking the same countable intersections as in the proof of \eqref{eq-squarefree-quantitative} proves \eqref{eq-squarefree-two-powers}.
    \end{proof}

    \section*{Acknowledgements}
    The author thanks Imre Z. Ruzsa for helpful discussions on the representations of almost all odd integers as the sum of a squarefree integer and a power of two. The author also acknowledges the use of ChatGPT in the preparation of this manuscript.

\end{document}